\input amstex\documentstyle{amsppt}  
\pagewidth{12.5cm}\pageheight{19cm}\magnification\magstep1
\topmatter
\title The Grothendieck group of unipotent representations: a new basis\endtitle
\author G. Lusztig\endauthor
\address{Department of Mathematics, M.I.T., Cambridge, MA 02139}\endaddress
\thanks{Supported by NSF grant DMS-1855773.}\endthanks
\endtopmatter   
\document

\define\bx{\boxed}

\define\Irr{\text{\rm Irr}}

\define\mpb{\medpagebreak}

\define\frl{\forall}
\define\pe{\perp}
\define\si{\sim}
\define\wt{\widetilde}
\define\sqc{\sqcup}

\define\qua{\quad}

\define\op{\oplus}
   
\define\part{\partial}
\define\emp{\emptyset}

\define\ra{\rangle}
\define\n{\notin}

\define\m{\mapsto}
\define\do{\dots}
\define\la{\langle}

\define\sub{\subset}    
\define\bxt{\boxtimes}
\define\T{\times}
\define\ti{\tilde}
\define\nl{\newline}
\redefine\i{^{-1}}

\define\un{\underline}

\define\ot{\otimes}

\define\Ind{\text{\rm Ind}}

\define\tr{\text{\rm tr}}

\define\di{\diamond}
\redefine\spa{\spadesuit}

\define\a{\alpha}

\redefine\d{\delta}
\define\e{\epsilon}

\define\p{\pi}

\define\ps{\psi}
\define\r{\rho}
\define\s{\sigma}
\redefine\t{\tau}
\define\th{\theta}
\define\k{\kappa}
\redefine\l{\lambda}
\define\z{\zeta}
\define\x{\xi}

\redefine\G{\Gamma}

\redefine\L{\Lambda}

\define\Ps{\Psi}

\redefine\ss{\bold s}

\define\BB{\bold B}
\define\CC{\bold C}

\define\FF{\bold F}

\define\NN{\bold N}

\define\RR{\bold R}
\define\SS{\bold S}

\define\ZZ{\bold Z}

\define\cf{\Cal F}
\define\cg{\Cal G}

\define\ci{\Cal I}

\define\cu{\Cal U}

\define\cx{\Cal X}

\define\tx{\ti x}

\define\tB{\ti B}

\define\tF{\ti F}

\define\tI{\ti I}

\define\EEIGHT{L1}
\define\SYMP{L2}
\define\ORA{L3}
\define\LEA{L4}
\define\NEW{L5}
\head Introduction\endhead
\subhead 0.1\endsubhead
Let $G$ be an adjoint simple algebraic group defined and split over a finite field $\FF_q$ and let $G(\FF_q)$ be
the finite group of $\FF_q$-rational points of $G$. Let $W$ be the Weyl group of $G$. 
We fix a family $c$ (in the sense of \cite{\EEIGHT})  in the set of irreducible representations of $W$. 
(This is the same as fixing a two-sided cell of $W$.) To $c$ we associate a finite group $\cg_c$ and an 
imbedding $c\sub M(\cg_c)$  (with image $M_0(\cg_c)$) as in \cite{\EEIGHT}, \cite{\ORA}. 
Here for any finite group $\G$, $M(\G)$ consists of pairs $(x,\r)$ where $x\in\G$ and $\r$ is an irreducible 
representation of the centralizer of $x$; these pairs are taken up to $\G$-conjugacy; let
$\CC[M(\G)]$ be the $\CC$-vector space with basis $M(\G)$ and let $A_\G:\CC[M(\G)]@>>>\CC[M(\G)]$ be the
``non-abelian Fourier transform'' (as in \cite{\EEIGHT}).
An element $f\in\CC[M(\G)]$ is said to be $\ge0$ if $f$ is a linear combinations of basis elements
$(x,\r)\in M(\G)$ with all coefficients in $\RR_{\ge0}$. As in \cite{\NEW} we say that $f\in\CC[M(\G)]$ is 
{\it bipositive} if $f\ge0$ and $A_\G(f)\ge0$.
 
Taking $\G=\cg_c$, we denote by $\CC[M_0(\cg_c)]$ the subspace of $\CC[M(\cg_c)]$ spanned by $M_0(\cg_c)$. 
In this paper we construct a new basis $\ti\BB_c$ of $\CC[M(\cg_c)]$. Here are some of the properties
of $\ti\BB_c$.

(I) All elements of $\ti\BB_c$ are bipositive.

(II) There is a unique bijection $M(\cg_c)@>\si>>\ti\BB_c$, $(x,\r)\m\widehat{(x,\r)}$, 
such that any 
$(x,\r)$ appears with 
nonzero coefficient in $\widehat{(x,\r)}$ ; this coefficient is actually $1$.

(III) Let $\le$ be transitive relation on $M(\cg_c)$ generated by the relation
for which $(x,\r),(x',\r')$ are related if 
$(x,\r)$ appears with nonzero coefficient in $\widehat{(x',\r')}$. Then $\le$ is a 
partial order on $M(\cg_c)$ in which $(1,1)$ is the unique minimal element.
In particular, the basis $\ti\BB_c$ is related to the basis $M(\cg_c)$ of 
$\CC[M(\cg_c)]$ by an upper triangular matrix with $1$ on diagonal and with integer 
entries.

(IV) $(1,1)$ appears with coefficient $1$ in any element of $\ti\BB_c$.

(V) The intersection $\ti\BB_c\cap\CC[M_0(\cg_c)]$ is the basis $\BB_c$ of $\CC[M_0(\cg_c)]$ defined in 
\cite{\NEW}.
\nl
Note that (III) and (V) imply 0.5(i) of \cite{\NEW} which was stated there without proof.

\subhead 0.2\endsubhead
Let $H\sub H'$ be subgroups of $\cg_c$ with $H$ normal in $H'$.
In 3.1 we define a linear map $\ss_{H,H'}:\CC[M(H'/H)]@>>>\CC[M(\G)]$ which commutes with the
non-abelian Fourier transform and takes bipositive elements to bipositive elements.

In the case where $G$ is of exceptional type,
our basis $\ti\BB_c$ is obtained by applying $\ss_{H,H'}$ to a very restricted set of bipositive elements 
(said to be {\it primitive}) of $\CC[M(H'/H)]$ where $H,H'$ are in the set of 
subgroups of $\cg_c$ which are either $\{1\}$ or are
associated in \cite{\LEA} to the various left cells of $W$ corresponding to $c$. 
This generalizes the definition of $\BB_c$ given in \cite{\NEW} where the linear map 
$\ss_{H,H'}$ was applied only to $(1,1)$.
In this case, our results can be interpreted as giving a new parametrization of $M(\cg_c)$
by triples $(H,H',\Xi)$ where $H,H'$ are as above and $\Xi$ runs through the primitive
bipositive elements of $\CC[M(H'/H)]$.
(In each case $H'/H$ is a symmetric group of small order.)

In the case where $G$ is of classical type our basis $\ti\BB_c$ will be defined using a
somewhat different approach. We will show elsewhere (based on results in \cite{\NEW,\S2}) that the approach
described above for exceptional types works also for classical types, leading to the same $\ti\BB_c$.

\subhead 0.3\endsubhead
Let $\Irr_c$ be the set of isomorphism classes of irreducible complex representations of $G(\FF_q)$ which are 
unipotent and are associated to $c$ as in \cite{\ORA}. Let $\cu_c$ be the (abelian) category of finite 
dimensional complex representations of $G(\FF_q)$ which are direct sums of representations in $\Irr_c$ and 
let $K_c$ be the Grothendieck group of $\cu_c$. In \cite{\ORA}, a bijection $M(\cg_c)@>\si>>\Irr_c$ is 
established. Via this bijection we can identify $\CC\ot K_c=\CC[M(\cg_c)]$ so that the basis $\Irr_c$
of $K_c$ becomes the basis $M[\cg_c]$ of $\CC[M(\cg_c)]$. Then the new basis $\ti\BB_c$ of $\CC[M(\cg_c)]$ 
becomes a new basis of $\CC\ot K_c$ (it also a $\ZZ$-basis of $K_c$).
The elements in this new $\ZZ$-basis of $K_c$ represent objects of $\cu_c$ which are called the
new (unipotent) representations of $G(\FF_q)$. They are in bijection with $\Irr_c$.
Note that taking disjoint union over the various families of $W$ we obtain a new basis for the Grothendieck
group of unipotent representations of $G(\FF_q)$.

In type $A_n$ we have $|c|=1$ and we can take $\ti\BB_c$ to consist of $(1,1)$; then the desired properties
of $\ti\BB_c$ are trivial. The properties above of $\ti\BB_c$ are verified in type $B_n,C_n,D_n$ in 
\S1. Another approach in type $D_n$ is sketched in \S2. The exceptional types are considered in \S3. 

\subhead 0.4\endsubhead
{\it Notation.} For $a,b$ in $\ZZ$ we set $[a,b]=\{z\in\ZZ;a\le z\le b\}$. 
For $a,b$ in $\ZZ$ we write
$a=_2b$ instead of $a=b\mod2$ and $a\ne_2b$ instead of $a\ne b\mod2$.
For a finite set $Y$ let $|Y|$ be the cardinal of $Y$.

\head 1. The set $\SS_D$\endhead
\subhead 1.1\endsubhead
Let $D\in\NN$. A subset $I$ of $[1,D]$ is said to be an {\it interval} if $I=[a,b]$ for some 
$a\le b$ in $[1,D]$. Let $\ci_D$ be the set of intervals of 
$[1,D]$. 
For $I=[a,b],I'=[a',b']$ in $\ci_D$ we write $I\prec I'$ whenever $a'<a\le b<b'$. 
We say that $I,I'$ are non-touching (and we write $I\spa I'$) if $a'-b\ge2$ or $a-b'\ge2$.
Let $R_D$ be the set
whose elements are the subsets of $\ci_D$. Let $\emp\in R_D$ be the empty subset of 
$\ci_D$.
For $B\in R_D$ and $h\in\{0,1\}$ we set $B^h=\{I\in B;|I|=_2h\}$. 

For $B\in R_D$ and $[a,b]\in\ci_D$ we define $\cx_B[a,b]=\cup_{I\in B^1;I\sub[a,b]}I$.

Let $I\in\ci_D$. A subset $E$ of $I$ is said to be {\it discrete} if $i\ne j$ in $E$ implies 
$i-j\ne\pm1$. Such $E$ is said to be maximal if $|E|=|I|/2$ (with $|I|$ even) or $|E|=(|I|+1)/2$ 
(with $|I|$ is odd). A maximal discrete subset of $I$ exists; it is unique if $|I|$ is odd. 

When $D\ge2$ and $i\in[1,D]$ we define an (injective) map $\x_i:\ci_{D-2}@>>>\ci_D$ by
$$\align&\x_i([a',b'])=[a'+2,b'+2]\text{ if }i\le a',\qua
\x_i([a',b'])=[a',b']\text{ if }i\ge b'+2, \\&
\x_i([a',b'])=[a',b'+2]\text{ if }a'<i<b'+2.\tag a\endalign$$ 
We define $t_i:R_{D-2}@>>>R_D$ by $B'\m\{\x_i(I');I'\in B'\}\sqc\{i\}$. We have $|t_i(B')|=|B'|+1$.

\subhead 1.2\endsubhead
We now assume that $D$ is even.
We say that $B\in R_D$ is {\it primitive } if it is of the form   

(a) $B=\{[1,D],[2,D-1],\do,[k,D+1-k]\}$ for some $k\in\NN$, $k\le D/2$.
\nl
For example, $B=\emp\in R_D$ is primitive (with $k=0$).
We define a subset $\SS_D$ of $R_D$ by induction on $D$ as follows.

If $D=0$, $\SS_D$ consists of a single element namely $\emp\in R_D$.
If $D\ge2$ we say that $B\in R_D$ is in $\SS_D$ if either $B$ is primitive, or

(b) there exists $i\in[1,D]$ and $B'\in\SS_{D-2}$ such that $B=t_i(B')$.
\nl
(This generalizes the definition of the set $S_D$ in \cite{\NEW,1.2} which can be viewed
as a subset of $\SS_D$.)

Let $\t_D:[1,D]@>>>[1,D]$ be the involution $i\m D+1-i$. It induces an involution $I\m\t_D(I)$ of 
$\ci_D$.
One can verify that $I\m\t_D(I)$ defines an involution $\SS_D@>>>\SS_D$; we denote it again by 
$\t_D$.

\subhead 1.3\endsubhead
For $D\ge0$, let $\SS_D^{prim}=\{B\in\SS_D;B\text{ primitive}\}$. 

Let $B\in R_D$. We consider the following properties $(P_0),(P_1),(P_2)$ that $B$ may or may not have.

$(P_0)$ {\it If $I\in B$, $\tI\in B$, then either $I=\tI$, or $I\spa\tI$, or $I\prec\tI$, or 
$\tI\prec I$.}

$(P_1)$ {\it If $[a,b]\in B^1$ and $b-a\ge2$ then $\cx_B[a+1,b-1]$ contains the unique maximal
discrete subset of $[a+1,b-1]$, that is, $\{a+1,a+3,a+5,\do,b-1\}$.}

$(P_2)$ {\it Let $k=|B^0|\in\NN$. There exists a (necessarily unique) sequence of integers 
$0=h_0<h_1<h_2<\do<h_{2k}<h_{2k+1}=D+1$ such that $B^0$ consists of $[h_1,h_{2k}]$, 
$[h_2,h_{2k-1}]$, $\do$, $[h_k,h_{k+1}]$. We have $h_j=_2j$ for $j\in[0,2k+1]$. Assume now that 
$k\ge1$ and that $j\in[0,2k]-\{k\}$ satisfies $h_{j+1}\ge h_j+3$. If $j\in[0,k-1]$, then 
$\cx_B[h_j+1,h_{j+1}-2]$ contains the unique maximal discrete subset of $[h_j+1,h_{j+1}-2]$; if 
$j\in[k+1,2k]$ then $\cx_B[h_j+2,h_{j+1}-1]$ contains the unique maximal discrete subset of 
$[h_j+2,h_{j+1}-1]$.}
\nl
Assume now that $D\ge2$, $i\in[1,D]$, $B'\in R_{D-2}$, $B=t_i(B')\in R_D$.
From the definitions we see that the following holds.

(a) {\it $B'$ satisfies $(P_0),(P_1),(P_2)$ if and only if $B$ satisfies $(P_0),(P_1),(P_2)$.}
\nl
Let $\SS'_D$ be the set of all $B\in R_D$ which satisfy $(P_0),(P_1),(P_2)$. 
(This generalizes the definition of the set $S'_D$ in \cite{\NEW, 1.3}.
Properties like $(P_0),(P_1)$ appeared in \cite{\NEW, 1.3}.)

In the setup of (a) we have the following consequence of (a).

(b) {\it We have $B'\in\SS'_{D-2}$ if and only if $B\in\SS'_D$.}
\nl
We show (extending \cite{\NEW, 1.3(c)}:

(c) {\it We have $\SS_D=\SS'_D$. In particular any $B\in\SS_D$ satisfies $(P_0),(P_1),(P_2)$.}
\nl
We argue by induction on $D$. If $D=0$, $\SS'_D$ consists of the empty set hence (c) holds in this 
case. Assume now that $D\ge2$. Let $B\in\SS_D$. We show that $\BB\in\SS'_D$. If $B\in\SS_D^{prim}$
then $B$ clearly is in $\SS'_D$. If $B\n\SS_D^{prim}$ then $B=t_i(B')$ for some 
$i,B'\in S_{D-2}$ as in 1.2(b). By the induction hypothesis we have $B'\in\SS'_{D-2}$. By (b) we have
$B\in\SS'_D$. We see that $B\in\SS_D\implies B\in\SS'_D$. Conversely, let $B\in\SS'_D$. We show that 
$B\in\SS_D$. If $B\in\SS_D^{prim}$ this is obvious. Thus we can assume that $B\n\SS_D^{prim}$.
From $(P_2)$ we see that $B^1\ne\emp$. Let $[a,b]\in B^1$ be such that $b-a$ is minimum. If $a<z<b$,  
$z=_2a+1$ then by $(P_1)$ we have $z\in[a',b']$ with $[a',b']\in B^1$, $b'-a'<b-a$, contradicting 
the minimality of $b-a$. We see that no $z$ as above exists. Thus, $[a,b]=\{i\}$ for some $i\in[1,D]$. 
Using $(P_0)$ and $\{i\}\in B$, we see that $B$ does not contain any interval of the form $[a,i]$ with 
$[a,i]$ with $a<i$, or $[i,b]$ with $i<b$, or $[a,i-1]$ with $a<i$ or $[i+1,b]$ with $i<b$; hence any 
interval of $B$ other than $\{i\}$ is of the form $\x_i[a',b']$ where $[a',b']\in\ci_{D-2}$. Thus we 
have $B=t_i(B')$ for some $B'\in R_{D-2}$. From (a) we deduce that $B'\in\SS'_{D-2}$. Using the 
induction hypothesis we deduce that
$B'\in\SS_{D-2}$. By the definition of $\SS_D$, we have $B\in\SS_D$. This completes the proof of (c).

\mpb

We show:

(d) {\it Let $B\in \SS_D$. If $I\in B^1,J\in B^0$ then $J\not\sub I$.}
\nl
We argue by induction on $|I|$. Let $I=[a,b]$, $J=[a',b']$. Assume that $J\sub I$. By $(P_0)$ we 
have $J\prec I$. Since $b'-a'$ is odd, then either $x=a'$ or $x=b'$ satisfies $x=_2a+1$. By $(P_1)$ 
we can find $I'\in B^1$ such that $I'\prec I$, $x\in I'$. We have $|I'|<|I|$. By the induction 
hypothesis we have $J\not\sub I'$. We have $I'\cap J\ne\emp$ and $I'\ne J$ hence $I'\prec J$ so that 
$x\n I'$, a contradiction. This proves (d).

\mpb

Let $B\in\SS_D$ and let $h_0<h_1<\do<h_{2k+1}$ be attached to $B$ as in $(P_2)$. We show:

(e) {\it If $[a,b]\in B^1$, then for some $j\in[0,2k]$ we have $h_j<a\le b<h_{j+1}$.}
\nl
We can find $j\in[0,2k]$ such that $h_j\le a\le h_{j+1}$. Assume first that $j\in[0,k-1]$.
Then $[h_j,h_{2k+1-j}]\cap[a,b]\ne\emp$ and $[h_j,h_{2k+1-j}]\ne[a,b]$ (one is in $B^0$, the other in 
$B^1$). Using $(P_0)$, we deduce $[h_j,h_{2k+1-j}]\prec[a,b]$ (which contradicts (d)) or 
$[a,b]\prec[h_j,h_{2k+1-j}]$ so that $h_j<a$. If $b\ge h_{2k-j}$, then $[h_{j+1},h_{2k-j}]\sub[a,b]$ 
contradicting (d). Thus we have $b<h_{2k-j}$. If $b\ge h_{j+1}$, then 
$[h_{j+1},h_{2k-j}]\cap[a,b]\ne\emp$ (it contains $b$) and $[h_{j+1},h_{2k-j}]\ne[a,b]$. Hence, by 
$(P_0)$, we have either $[h_{j+1},h_{2k-j}]\prec[a,b]$ (which again contradicts (d)) or 
$[a,b]\prec[h_{j+1},h_{2k-j}]$ hence $a>h_{j+1}$, contradicting our assumption. We see that 
$b<h_{j+1}$.

Assume next that $j\in[k+1,2k]$. Then $[h_{2k-j},h_{j+1}]\cap[a,b]\ne\emp$ (it contains $a$) and 
$[h_{2k-j},h_{j+1}]\ne[a,b]$ (one is in $B^0$, the other in $B^1$). Using $(P_0)$, we deduce 
$[h_{2k-j},h_{j+1}]\prec[a,b]$ (which contradicts (d)) or $[a,b]\prec[h_{2k-j},h_{j+1}]$, so that 
$b<h_{j+1}$. If $a=h_j$ then $[h_{2k+1-j},h_j]\cap[a,b]\ne\emp$ (it contains $a$) and 
$[h_{2k+1-j},h_j]\ne[a,b]\ne\emp$. Using $(P_0)$ we deduce $[h_{2k+1-j},h_j]\prec[a,b]$ (which 
contradicts (d)) or $[a,b]\prec[h_{2k+1-j},h_j]$ hence $a<h_j$, a contradiction. We see that $a>h_j$.

Finally, we assume that $j=k$. Then $[h_k,h_{k+1}]\cap[a,b]\ne\emp$ and $[h_k,h_{k+1}]\ne[a,b]$ (one is 
in $B^0$, the other in $B^1$). Using $(P_0)$, we deduce $[h_k,h_{k+1}]\prec[a,b]$ (which contradicts 
(d)) or $[a,b]\prec[h_k,h_{k+1}]$, so that $h_k<a$ and $b<h_{k+1}$. This proves (e).

The following result has already been proved as a part of the proof of (c).

(f) {\it Assume that $D\ge2$, $i\in[1,D]$. Let $B\in\SS_D$ be such that $\{i\}\in B$. 
Then there exists $B'\in\SS_{D-2}$ such that $B=t_i(B')$.}

\mpb

Let $B\in\SS_D$ and let $I=[a,b]\in B^1$. Let $\cx(I)=\{I'\in B^1;I'\sub I\}$. We show:

(g) {\it $|\cx(I)|=(b-a+2)/2$.}
\nl
We argue by induction on $|I|$. If $|I|=1$ then $\cx(I)=\{I\}$ and the result is clear.
Assume now that $|I|\ge3$. 
By $(P_0),(P_1)$ we can find $a=z_0<z_1<\do<z_r=b$ ($r\ge0$) such that
$z_0,z_1,\do, z_r$ are all congruent to $a\mod2$ and
$[z_0+1,z_1-1]\in B^1,[z_1+1,z_2-1]\in B^1,\do,[z_{r-1}+1,z_r-1]\in B^1$;
moreover, any $I'\in B^1$ such that $I'\prec I$ is contained in exactly one of
$[z_0+1,z_1-1],[z_1+1,z_2-1],\do,[z_{r-1}+1,z_r-1]$. It follows that
$|\cx)I)|=1+\sum_{j\in[0,r-1]}|\cx([z_j+1,z_{j+1}-1])|$. Using the induction hypothesis we can rewrite
the last equality as
$|\cx(I)|=1+\sum_{j\in[0,r-1]}((z_{j+1}-1)-(z_j+1)+2)/2=1+(b-a)/2$.
This proves (g).

\subhead 1.4\endsubhead
For $B\in\SS_D$, $h\in\{0,1\}$, $j\in[1,D]$ we set $B^h_j=\{I\in B^h;j\in I\}$. From the definitions
we deduce:

(a) {\it Assume that $D\ge2$, $i\in[1,D]$ and that $B'\in\SS_{D-2}$. Let $B=t_i(B')\in\SS_D$. Then 
$|B^0|=|B'{}^0|$. Moreover, for $h\in\{0,1\}$ and $r\in[1,D-2]$ we have:

$|B'_r{}^h)|=|B_r^h|$ if $r\le i-2$, $|B'_r{}^h|=|B_{r+2}^h|$ if $r\ge i$,

$|B^h_{i-1}|=|B^h_{i+1}|=|B'_{i-1}{}^h|$, $|B^h_i|=|B'_{i-1}{}^h|+h$ if $1<i<D$,

$|B^h_{i-1}|=0$ if $i=D$, $|B^h_{i+1}|=0$ if $i=1$.}

This extends \cite{\NEW, 1.4(a)}.

\subhead 1.5\endsubhead
Let $B\in\SS_D-\SS_S^{prim}$. 
As we noted in the proof of 1.3(c), in this case we must have $B^1\ne\emp$ and we have $\{j\}\in B^1$
for some $j\in[1,D]$; we assume that $j$ is as small as possible (then it is uniquely determined).
As in that proof we have $B=t_j(B')$ where $B'\in\SS_{D-2}$. Let $i$ be the smallest number in 
$\cup_{I\in B^1}I$. We have $i\le j$. We show:

(a) {\it For any $h\in[i,j]$, we have $[h,\ti h]\in B^1$ for a unique $\ti h\in[h,D]$; moreover we 
have $j\le\ti h$.}
\nl
We argue by induction on $D$. When $D\le1$ the result is obvious. We now assume that $D\ge2$.
Assume first that $i=j$. By $(P_0)$, $\{j\}\in B^1$ implies that we cannot have $[j,b]\in B^1$ with
 $j<b$; thus (a) holds in this case. We can assume that $i<j$. We have $[i,b]\in B^1$ for some $b>i$ 
hence $|B^1|\ge2$ so that $|B'{}^1|\ge1$ and $B'\n\SS_{D-2}^{prim}$. Then $i',j'$ are defined in 
terms of $B'$ in the same way as $i,j$ are defined in terms of $B$. From $(P_1)$ we see that there 
exists $j_1$ such that $i<j_1<b$ and such that $\{j_1\}\in B$. By the minimality of $j$ we must have 
$j\le j_1$. Thus we have $i<j<b$. We have $[i,b]=\x_j[i,b-2]]$ hence $[i,b-2]\in B'{}^1$. This 
implies that $i'\le i$. We have $[i',c]\in B'{}^1$ for some $c\in[i',D-2]$, $c=_2i'$; hence 
$[i',c']\in B^1$ for some $c'\ge i'$ so that $i'\ge i$. Thus we have $i'=i$. By the induction 
hypothesis, the following holds:

(b) {\it For any $r\in[i,j']$, we have $[r,r_1]\in B'{}^1$ for a unique $r_1$; moreover $j'\le r_1$.}
\nl
If $j'\le j-2$ then $\{j'\}=\x_j(\{j'\})\in B$. Hence $j'\ge j$ by the minimality of $j$; this is a 
contradiction. Thus we have $j'\ge j-1$.

Let $r\in[i,j-1]$. Then we have also $r\in[i,j']$ hence $r_1$ is defined as in (b). We have
$[r,r_1]\in B'{}^1$ hence $[r,r_1+2]\in B^1$ (we use that $r<j\le j'+1\le r_1+1<r_1+2$); we have 
$j<r_1+2$. Assume now that $[r,r_2]\in B^1$ with $r\le r_2$. Then $r<r_2$ (by the minimality of $j$).
If $j=r_2$ or $j=r_2+1$ then applying $(P_0)$ to $\{j\},[r,r_2]$ gives a contradiction.
Thus we must have either $r<j<r_2$ or $j>r_2+1$. If $j>r_2+1$ then $[r,r_2]\in B'{}^1$ hence by (b), 
$r_2=r_1$, hence $j>r_1+1$ contradicting $j<r_1+2$. Thus we have $r<j<r_2$, so that 
$[r,r_2-2]\in B'{}^1$ hence by (b), $r_2-2=r_1$. Thus we have $r<j<r_2$ so that $[r,r_2-2]\in B'{}^1$ 
hence by (b), $r_2-2=r_1$.

Next we assume that $r=j$. In this case we have $\{r\}\in B^1$. Moreover, if $[r,r']\in B^1$ with
$r\le r'\le D$, then we cannot have $r<r'$ (if $r<r'$ then applying $(P_0)$ to $\{r\},[r,r']$ gives a
contradiction). This proves (a).

We show:

(c) {\it Assume that $j<D$ and that $i\le h<j$. Then $\ti h$ in (a) satisfies $\ti h>j$.}
\nl
Assume that $\ti h=j$, so that $[h,j]\in B^1$. Since $h<j$, 
applying $(P_0)$ to $\{j\},[h,j]$ gives a contradiction.  This proves (c).

We show:

(d) {\it Assume that $j<D$ and that $r\in[j+1,D]$. We have $[j+1,r]\n B^1$.}
\nl
Assume that $[j+1,r]\in B^1$. Applying $(P_0)$ to $\{j\},[j+1,r]$ gives a contradiction.  This proves 
(d).

We show:

(e) {\it For $h\in[i,j]$ we have $|B^1_h|=h-i+1$. If $j<D$ we have $|B^1_{j+1}|=j-i$.}
\nl
Let $h\in[i,j]$. Then for any $h'\in[i,h]$, $B^1_h$ contains $[h',\ti h']$ (since $h\le\ti h'$), see 
(a). Conversely, assume that $[a,b]\in B^1_h$. We have $a\le h$. 
By the definition of $i$ we have $i\le a$.
By the uniqueness statement in (a) we have $b=\ti a$ so that $[a,b]$ is one of the $h-i+1$ intervals 
$[h',\ti h']$ above. This proves the first assertion of (e). Assume now that $j<D$. If $h'\in[i,j]$, 
$h'<j$,  then $[h',\ti h']\in B^1_{j+1}$, by (c). Conversely, assume that
$[a,b]\in B^1_{j+1}$. We have $a\le j+1$ and by (d) we have $a\ne j+1$ so that $a\le j$. If $a=j$, then 
by the uniqueness in (a) we have $b=j$ which contradicts $j+1\in[a,b]$. Thus we have $a\le j-1$. We 
see that $[a,b]$ is one of the $j-i$ intervals $[h',\ti h']$ with $h'\in[i,j]$, $h'<j$. This proves (e).

We show:

(f) {\it Let $e=|B^0_i|$. For $h\in[i,j]$ we have $|B^0_h|=e$. If $j<D$ we have $|B^0_{j+1}|=e$.}
\nl
Let $I\in B^0_i$. Since $I$ and $[i,\ti i]$ are not disjoint and not equal, we must have 
$[i,\ti i]\prec I$ or $I\prec[i,\ti i]$ (this last case cannot occur since $i\in[i,\ti i]$). Thus we 
have $[i,\ti i]\prec I$. We have $[i,j]\sub[i,\ti i]$ hence $[i,j]\prec I$ so that $I\in B^0_h$ for 
any $h\in[i,j]$. If in addition $j<D$, then from $[i,j]\prec I$ we deduce $[i,j+1]\sub I$ so that 
$I\in B^0_{j+1}$. Conversely, assume that $h\in[i,j]$ and $I'\in B^0_h$. Since $I'$ and $[h,\ti h]$ are 
not disjoint and not equal, we must have $[h,\ti h]\prec I'$ or $I'\prec[h,\ti h]$ (this last case 
cannot 
occur since $h\in I'$). Thus we have $[h,\ti h]\prec I'$. If $i<h$, it follows that $h-1\in I'$ so that 
$[h-1,\wt{h-1}]\prec I'$. Repeating this argument we see that $[h'',\ti h'']\prec I'$ for any 
$h''\in[i,h]$, so that in particular we have $i\in I'$ and $I'\in B^0_i$. If in addition $j<D$ and 
$I'\in B^0_{j+1}$, then $I',[j,\ti j]$ are not non-touching and are not equal hence we must have 
$[j,\ti j]\prec I'$ or $I'\prec[j,\ti j]$ (this last case cannot occur since it contradicts
1.3(d)). Thus we have $[j,\ti j]\prec I'$ which by 
the earlier part of the proof implies $I'\in B^0_i$. This proves (f).

\subhead 1.6\endsubhead
For any $n\in\NN$ we define $\un n\in\{0,1\}$ by $n=_2\un{n}$.
For $B\in\SS_D$, $j\in[1,D]$, we set $\k=\un{|B^0|}$ and  
$$f_j(B)=|B^1_j|-|B^0_j|-\k\in\ZZ,$$
$$\e_j(B)=f_j(B)(f_j(B)+1)/2\in\FF_2.$$
This extends a definition in \cite{\NEW, 1.6}. We have 

$\e_j(B)=1$ if $f_j(B)\in(4\ZZ+1)\cup(4\ZZ+2)$, $\e_j(B)=0$ if $f_j(B)\in(4\ZZ+3)\cup(4\ZZ)$. 

Assume now that $B\n S_D^{prim}$. Let $i\le j$ in $[1,D]$ be as in 1.5. Let $e=|B^0_i|+\k$.
From 1.5(e),(f) we deduce:

(a) {\it We have 
$$(f_i(B),f_{i+1}(B),\do,f_j(B))=(1-e,2-e,3-e,\do,j-i-e,j-i+1-e).$$ 
If $j<D$, we have $f_{j+1}(B)=j-i-e$.}
\nl
From (a) we deduce:
$$\align&(\e_i(B),e_{i+1}(B),\do,\e_j(B))=((1-e)(2-e)/2,(2-e)(3-e)/2,\\&
(3-e)(4-e)/2,\do,(j-i-e)(j-i-e+1)/2,(j-i-e+1)(j-i-e+2));\tag b\endalign$$
(c) {\it if $j<D$, then} $\e_{j+1}(B)=(j-i-e)(j-i-e+1)/2)$.
\nl
This extends \cite{\NEW, 1.6(b),(c)}.

For future reference we note:

(d) {\it If $c\in\ZZ$ then $c(c+1)/2\ne_2(c+2)(c+3)/2$.}

(e) {\it If $c\in2\ZZ$ then $c(c+1)/2\ne_2(c+1)(c+2)/2$.}

\subhead 1.7\endsubhead
Let $B\in\SS_D$, $\tB\in\SS_D$ be such that $B,\tB$ are not primitive and $\e_h(B)=\e_h(\tB)$
for any $h\in[1,D]$ and $|B^0|=|\tB^0|$. We show (extending \cite{\NEW, 1.7(a)}:

(a) {\it We can find $z\in[1,D]$ such that $\{z\}\in B$, $\{z\}\in\tB$.}
\nl
Let $h_0<h_1<h_2<\do<h_{2k}<h_{2k+1}$ be the sequence attached to $B$ in $(P_2)$; let
$\ti h_0<\ti h_1<\ti h_2<\do<\ti h_{2\ti k}<\ti h_{2\ti k+1}$ 
be the analogous sequence attached to $\tB$. Here $k=\ti k=|B^0|=|\tB^0|$. We shall need the following 
preparatory result.

(b) {\it Assume that $s\in[0,k-1]$ is such that 
$$(h_0,h_1,\do,h_s)=(\ti h_0,\ti h_1,\do,\ti h_s)=(0,1,\do,s).$$
Then either $h_{s+1}=\ti h_{s+1}=s+1$, or the conclusion of (a) holds.}
\nl
Let $i\le j$ be attached to $B$ as in 1.5. Let $\ti i\le\ti j$ be similarly attached to $\tB$.
Assume first that $h_{s+1}>s+1$, $\ti h_{s+1}=s+1$. We have $|B^0_{s+1}|=s$, $|\tB^0_{s+1}|=s+1$, 
$|\tB^1_{s+1}|=0$ (we use 1.3(e)) and by $(P_2)$ we have $|B^1_{s+1}|\ge1$. We see that $i=s+1$ and 
from 1.5(e) we have $|B^1_{s+1}|=1$. Thus, $f_{s+1}(B)=1-s-\k$, $f_{s+1}(\tB)=-1-s-\k$
(where $\k=\un k$), so 
that $\e_{s+1}(B)=(1-s-\k)(2-s-\k)/2$, $\e_{s+1}(\tB)=(-1-s-\k)(-s-\k)/2$. It 
follows that 
$(1-s-\k)(2-s-\k)/2=_2(-1-s-\k)(-s-\k)/2$, contradicting 1.6(d). Thus, if 
$\ti h_{s+1}=s+1$, then $h_{s+1}=s+1$. Similarly, if $h_{s+1}=s+1$, then $\ti h_{s+1}=s+1$. Assume now 
that $h_{s+1}>s+1$ and $\ti h_{s+1}>s+1$. By $(P_2)$ we have $i=\ti i=s+1$. If $j<\ti j$, then $j<D$ 
and from 1.6(b),(c), we see that 
$$\e_{j+1}(B)=(j-i-s-\k)(j-i-s-\k+1)/2, \e_{j+1}(\tB)=(j-i-s-\k+2)(j-i-s-\k+3)/2$$ 
so that
$$(j-i-s-\k)(j-i-s-\k+1)/2=_2(j-i-s-\k+2)(j-i-s-\k+3)/2,$$
contradicting 1.6(d). Thus we have $j\ge\ti j$. Similarly, we have $\ti j\ge j$. Hence $\ti j=j$, so 
that
(a) holds with $z=j=\ti j$. The only remaining case is that where $h_{s+1}=s+1$ and $\ti h_{s+1}=s+1$. 
This proves (b).

We shall need a second preparatory result.

(c) {\it Assume that $s\in[0,k-1]$ is such that 
$$(h_{2k-s+1},\do,h_{2k},h_{2k+1})=(\ti h_{2k-s+1},\do,\ti h_{2k},\ti h_{2k+1})=(D-s+1,\do,D,D+1).$$
Then either $h_{2k-s}=\ti h_{2k-s}=D-s$ or the conclusion of (a) holds.}
\nl
We note that the assumptions of (b) are satisfied when $B,\tB$ are replaced by 
$\t_D(B),\t_D(\tB)$ (see 1.2). Hence from (b) we deduce that either $h_{2k-s}=\ti h_{2k-s}=D-s$ or 
there exists $u\in[1,D]$ such that $\{u\}\in\t_D(B)$, $\{u\}\in\t_D(\tB)$ (which implies that
$\{\t_D(u)\}\in B$, $\{\t_D(u)\}\in\tB$. This proves (c). 

\mpb

Next we note that the assumption of (b) (and that of (c)) is satisfied when $s=0$. Hence from (b),(c) we 
obtain by induction on $s$ the following result.

(d) {\it We have either 
$$\align&
(h_0,h_1,\do,h_k,h_{k+1},\do,h_{2k+1})=(\ti h_0,\ti h_1,\do,\ti h_k,\ti h_{k+1},\do,\ti h_{2k+1})\\&
=(0,1,\do,k,D-k+1,\do,D,D+1)\endalign$$
or the conclusion of (a) holds.}
\nl  
Thus, to prove (a) we can assume that $B,\tB$ are as in the first alternative of (d). We have 
$k<i\le j<D-k+1$, $k<\ti i\le\ti j<D-k+1$ (we use 1.3(e) and $(P_1)$). 
Assume first that $j<\ti j$ (so that $j<D$) and $i<\ti i$.

 From 1.6 for $B$ we have 
$\e_i(B)=(1-k-\un k)(2-k-\un k)/2$. From $i<\ti i$ we have $\e_i(\tB)=(-k-\un k)(1-k-\un k)/2$. Thus
$(1-k-\un k)(2-k-\un k)/2=_2(-k-\un k)(1-k-\un k)/2$. This contradicts 1.6(e) since $k+\un k$ is even.
Thus we must have $i\ge\ti i$. 
Next we asssume that $j<\ti j$ (so that $j<D$) and $\ti i<i$. From 1.6 for $\tB$ we have 
$\e_{\ti i}(\tB)=(1-k-\un k)(2-k-\un k)/2$. From $\ti i<i$ we have 
$\e_{\ti i}(B)=(-k-\un k)(-k-\un k+1)/2$.
$(1-k-\un k)(2-k-\un k)/2=_2(-k-\un k)(1-k-\un k)/2$. This contradicts 1.6(e) since $k+\un k$ is even.
Thus, when $j<\ti j$ we must have $i=\ti i$. 
From 1.6(c) for $B$ we have $e_{j+1}(B)=(j-i-k-\un k)(j-i-k-\un k+1)/2$  and from 1.6(b) for $\tB$ we 
have $e_{j+1}(\tB)=(j-i-k-\un k+2)(j-i-k-\un k+3)/2$. It follows that 
$$(j-i-k-\un k)(j-i-k-\un k+1)/2)=_2(j-i-k-\un k+2)(j-i-k-\un k+3)/2,$$ 
contradicting 1.6(d). We see that $j<\ti j$ leads to a contradiction. Similarly, $\ti j<j$ leads to a 
contradiction. Thus we must have $j=\ti j$, so that (a) holds with $z=j=\ti j$. This completes the 
proof of (a).

\subhead 1.8\endsubhead
Let $B\in\SS_D$, $\tB\in\SS_D$.

(a) {\it Assume that $\tB\in\SS_D^{prim}$, that $\e_h(B)=\e_h(\tB)$ for any $h\in[1,D]$ and that
$|B^0|=|\tB^0|$. Then $\tB=B$.}
\nl
The proof is similar to that of 1.7(a). Assume that $B\n\SS_D^{prim}$. Let $i\le j$ be attached 
to $B$ as in 1.5. Let $h_0<h_1<h_2<\do<h_{2k}<h_{2k+1}$ be the sequence attached to $B$ in $(P_2)$; let
$\ti h_0<\ti h_1<\ti h_2<\do<\ti h_{2\ti k}<\ti h_{2\ti k+1}$ (that is, $0<1<\do<k<D+1-k<\do<D<D+1$) 
be the analogous sequence attached to $\tB$. We have $k=\ti k=|B^0|=|\tB^0|$.

\mpb

We show the following variant of 1.7(b).

(b) {\it Assume that $s\in[0,k-1]$ is such that $(h_0,h_1,\do,h_s)=(0,1,\do,s)$. Then $h_{s+1}=s+1$.}
\nl
Assume first that $h_{s+1}>s+1$. We have $|B^0_{s+1}|=s$, $|\tB^0_{s+1}|=s+1$, 
$|\tB^1_{s+1}|=0$ (we use 1.3(e)) and by $(P_2)$ we have $|B^1_{s+1}|\ge1$. We see that $i=s+1$ and from 
1.5(e) we have $|B^1_{s+1}|=1$. Thus, $f_{s+1}(B)=1-s-\un k$, $f_{s+1}(\tB)=-1-s-\un k$, so that 
$$\e_{s+1}(B)=(1-s-\un k)(2-s-\un k)/2, \e_{s+1}(\tB)=(-1-s-\un k)(-s-\un k)/2.$$
 It follows that $(1-s-\un k)(2-s-\un k)/2=_2(-1-s-\un k)(-s-\un k)/2$, contradicting 1.6(d).
Thus, we must have $h_{s+1}=s+1$. This proves (b).

Next we show the following variant of 1.7(c).

(c) {\it Assume that $s\in[0,k-1]$ is such that $(h_{2k-s+1},\do,h_{2k},h_{2k+1})=(D-s+1,\do,D,D+1)$. 
Then $h_{2k-s}=D-s$.}
\nl
We note that the assumptions of (b) are satisfied when $B,\tB$ are replaced by $\t_D(B),\t_D(\tB)$.
Hence from (b) we deduce that $h_{2k-s}=D-s$. This proves (c). 

Now we note that the assumption of (b) (and that of (c)) is satisfied when $s=0$. Hence from (b),(c) we 
obtain by induction on $s$ the following result.

(d) {\it We have }
$$\align&(h_0,h_1,\do,h_k,h_{k+1},\do,h_{2k+1})
=(\ti h_0,\ti h_1,\do,\ti h_k,\ti h_{k+1},\do,\ti h_{2k+1})\\&
=(0,1,\do,k,D-k+1,\do,D,D+1).\endalign$$
Using (d) and 1.6 we see that 
$e_i(B)=(1-k-\un k)(2-k-\un k)/2$. On the other hand we have $e_i(\tB)=(-k-\un k)(-k-\un k+1)/2$.
We get $(1-k-\un k)(2-k-\un k)/2=_2(-k-\un k)(-k-\un k+1)/2$, contradicting 1.6(e) since $k+\un k$ 
is even.
Thus $B\n\SS_D^{prim}$ leads to a contradiction. Thus both $B,\tB$ are primitive.
Since $B,\tB$ are primitive and $|B^0|=|\tB^0|$, we see that 
$B=\tB$. This proves (a).

\subhead 1.9\endsubhead
We no longer assume that $D$ is even. Let $V$ be the $\FF_2$-vector space with basis $\{e_i;i\in[1,D]\}$.
For any subset $I$ of $[1,D]$ let $e_I=\sum_{i\in I}e_i\in V$. We define a symplectic form 
$(,):V\T V@>>>\FF_2$ by 
$(e_i,e_j)=1$ if $i-j=\pm1$, $(e_i,e_j)=0$ if $i-j\ne\pm1$. This symplectic form is nondegenerate if $D$ is even
while if $D$ is odd it has a one dimensional radical spanned by $e_1+e_3+e_5+\do+e_D$.

For any subset $Z$ of $V$ we set $Z^\pe=\{x\in V;(x,z)=0\qua\frl z\in Z\}$.

When $D\ge2$ we denote by $V'$ the $\FF_2$-vector space with basis $\{e'_i;i\in[1,D-2]\}$.
For any $I'\sub[1,D-2]$ let $e'_{I'}=\sum_{i\in I'}e'_i\in V'$. 
We define a symplectic form $(,)':V'\T V'@>>>\FF_2$ by $(e'_i,e'_j)=1$ if $i-j=\pm1$,
$(e'_i,e'_j)=0$ if $i-j\ne\pm1$. 

When $D\ge2$, for any $i\in[1,D]$ there is a unique linear map $T_i:V'@>>>V$ such that 
the sequence $T_i(e'_1),T_i(e'_2),\do,T_i(e'_{D-2})$ is:

$e_1,e_2,\do,e_{i-2},e_{i-1}+e_i+e_{i+1},e_{i+2},e_{i+3},\do,e_D$ (if $1<i<D$),

$e_3,e_4,\do,e_D$ (if $i=1$),

$e_1,e_2,\do,e_{D-2}$ (if $i=D$).
\nl
Note that $T_i$ is injective and $(x,y)'=(T_i(x),T_i(y))$ for any $x,y$ in $V'$.
For any $I'\in\ci_{D-2}$ we have $T_i(e'_{I'})=e_{\x_i(I')}$. Let $V_i$ be the image of $T_i:V'@>>>V$. 
From the definitions we deduce:

(a) {\it $e_i^\pe=V_i\op \FF_2e_i$.}
\nl
{\it In the remainder of this section we assume that $D$ is even.}

If $D\ge2$, for $j\in[1,D-2]$ let $f'_j:\SS_{D-2}@>>>\ZZ$, $\e'_j:\SS_{D-2}@>>>\FF_2$ be the 
analogues of $f_i:\SS_D@>>>\ZZ$, $\e_i:\SS_D@>>>\FF_2$ when $D$ is replaced by $D-2$.

For $B\in\SS_D$, we define $\e(B)\in V$ by $\e(B)=\sum_{i\in[1,D]}\e_i(B)e_i$. 
If $D\ge2$, for $B'\in\SS_{D-2}$ we define $\e'(B')\in V'$ by 
$\e'(B')=\sum_{j\in[1,D-2]}\e'_j(B')e'_j$. We show (extending \cite{\NEW, 1.9(b)}:

(b) {\it Assume that $D\ge2$, $i\in[1,D]$. Let $B'\in\SS_{D-2},B=t_i(B')\in\SS_D$. Then
$\e(B)=T_i(\e'(B'))+ce_i$ for some $c\in\FF_2$.}
\nl
An equivalent statement is: for any $j\in[1,D]-\{i\}$ we have $\e_j(B)=\e'_{j'}(B')$ 
if $j'\in[1,D-2]$ 
is such that $j\in\x_i(\{j'\})$; and $\e_j(B)=0$ if no such $j'$ exists. It is enough to show:

$f'_h(B')=f_h(B)$ if $1\le h\le i-2$,

$f'_{h-2}(B')=f_h(B)$ if $i+2\le h\le D$,

$f_{i-1}(B)=f_{i+1}(B)=f'_{i-1}(B')$ if $1<i<D$,

$f_{i-1}(B)\in\{0,-1\}$ (hence $\e_{i-1}(B)=0$) if $i=D$, 

$f_{i+1}(B)\in\{0,-1\}$ (hence $\e_{i+1}(B)=0$) if $i=1$.
\nl
This follows from 1.4(a).

\mpb

For $B\in\SS_D$ let $\la B\ra$ be the subspace of $V$ generated by $\{e_I;I\in B\}$. 
For $B'\in\SS_{D-2}$ let $\la B'\ra$ be the subspace of $V'$ generated by $\{e'_{I'};I'\in B'\}$. We 
show (extending \cite{\NEW, 1.9(c)}:

(c) {\it Let $B\in\SS_D$. We have $\e(B)\in\la B\ra$.
If $D\ge2,i\in[1,D]$, $B'\in\SS_{D-2},B=t_i(B')\in\SS_D$, then $\la B\ra=T_i(\la B'\ra)\op\FF_2e_i$.}
\nl
To prove the first assertion of (c) we argue by induction on $D$. For 
$d=0$ there is nothing to prove. Assume that $d\ge1$. Let $i,B'$ be as 
in (b). By the induction hypothesis we have $\e'(B')\in\la B'\ra\sub V'$. Using (b) we see that it is 
enough to show that $T_i(\la B'\ra)\sub\la B\ra$. (Since $\{i\}\in B$, we have $e_i\in\la B\ra$.) 
Using the equality $T_i(e'_{I'})=e_{\x_i(I')}$ for any $I'\in B'$ it remains to note that $\x_i(I')\in B$ 
for $I'\in B'$. 
This proves the first assertion of (c). The same proof shows the second assertion of (c).

\mpb

For $s\in[0,D/2]$ we set $t=s/2$ if $s$ is even and $t=(s+1)/2$ if $s$ is odd; we denote by $V(s)$ 
the set of vectors $x\in V$ such that $x=e_{[a_1,b_1]}+e_{[a_2,b_2]}+\do+e_{[a_t,b_t]}$ with 
$[a_r,b_r]\in\ci_D$ with any two of them non-touching and with 
$a_1=_2a_2=_2\do=_2a_t=_2s$,  $b_1=_2b_2=_2\do=_2b_t=_2s+1$. 
For such $x$ we set $n(x)=a_1+b_1+a_2+ b_2+\do+a_t+b_t\in\NN$.

Assume for example that 
$$B=\{[1,D],[2,D-1],\do,[D/2,(D/2)+1]\}.$$
We have $|B^0_i|=i$ for $i\in[1,D/2]$, $|B^0_i|=D-i+1$ for $i\in[(D/2)+1,D]$,
$|B^1_i|=0$ for all $i$. It follows that 
$$\e(B)=e_{[2,3]}+e_{[6,7]}+e_{[10,11]}+\do+e_{[D-2,D-1]}\in V(D/2)\text{ if $D/2$ is even},\tag d$$
$$\e(B)=e_{[1,2]}+e_{[5,6]}+e_{[9,10]}+\do+e_{[D-1,D]}\in V(D/2)\text{ if $D/2$ is odd}.\tag e$$
More generally, assume that 
$$B=\{[1,D],[2,D-1],\do,[s,D+1-s]\}\text{ where }s\in[0,D/2].\tag f$$
We have $|B^0_i|=i$ for $i\in[1,s]$, $|B^0_i|=D-i+1$ for $i\in[D-s+1,D]$, $|B^1_i|=0$ for all $i$. It 
follows that 

If $s=0$ then $\e(B)=0$;
 
if $s=1$ then $\e(B)=e_{[1,D]}$; 

if $s=2$ then $\e(B)=e_{[2,D-1]}$;

if $s=3$ then $\e(B)=e_{[1,2]}+e_{[D-1,D]}$; 

if $s=4$ then $\e(B)=e_{[2,3]}+e_{[D-2,D-1]}$;

if $s=5$ then $\e(B)=e_{[1,2]}+e_{[5,D-4]}+e_{[D-1,D]}$; 

if $s=6$ then $\e(B)=e_{[2,3]}+e_{[6,D-5]}+e_{[D-2,D-1]}$, etc.
\nl
Thus,

(g) $\e(B)\in V(s)$.

\mpb

Let $B\in\SS_D$. Using $(P_0)$ we deduce:

(h) {\it $\la B\ra$ is an isotropic subspace of $V$.}
\nl
We show (extending \cite{\NEW, 2.1(b)}):

(i) {\it $\{e_I;I\in B\}$ is an $\FF_2$-basis of $\la B\ra$.}
\nl
We can assume that $D\ge2$. Assume that $\sum_{I\in B}c_Ie_I=0$ with $c_I\in\FF_2$ not all zero.
We can find $I_1=[a,b]\in B$ with $c_{I_1}\ne0$ and $|I_1|$ maximal. If $a\in I'$
with $I'\in B$, $I'\ne I_1$, $c_{I'}\ne0$, then by $(P_0)$ we have $I_1\prec I'$ (contradicting the
maximality of $|I_1|$) or $I'\prec I_1$ (contradicting $a\in I'$). Thus no $I'$ as above exists.
Thus when $\sum_{I\in B}c_Ie_I$ is written in the basis $\{e_j;j\in[1,2d]\}$, the coefficient of $e_a$
is $c_{I_1}$ hence $c_{I_1}=0$, contradicting $c_{I_1}\ne0$. Thus (i) holds for $B$.

\subhead 1.10\endsubhead
Let $B\in\SS_D,\tB\in\SS_D$. We show:

(a) {\it If $\e(B)=\e(\tB)$ and $|B^0|=|\tB^0|$, then $B=\tB$.}
\nl
Wec argue by induction on $D$. If $D=0$, there is nothing to prove. Assume that $D\ge2$. If 
$\tB\in\SS_D^{prim}$, then (a) follows from 1.8(a). Similarly, (a) holds if and $B\in\SS_D^{prim}$. 
Thus, we can assume that $B$ and $\tB$ are not primitive. By 1.7(a) we can find $i\in[1,D]$ 
such that 
$\{i\}\in B^1$, $\{i\}\in\tB^1$. By 1.3(f) we then have $B=t_i(B')$, $\tB=t_i(\tB')$ with 
$B'\in\SS_{D-2}$, $\tB'\in\SS_{D-2}$. Using our assumption and 1.9(b) we see that 
$T_i(\e'(B'))=T_i(\e'(\tB'))+ce_i$ for some $c\in\FF_2$. Using 1.9(a) we see that $c=0$ so that 
$T_i(\e'(B'))=T_i(\e'(\tB'))$. Since $T_i$ is injective, we deduce $\e'(B')=\e'(\tB')$.
We have also $|B'{}^0|=|\tB'{}^0|$. By the induction hypothesis we have $B'=\tB'$ hence $B=\tB$. This 
proves (a).

\subhead 1.11\endsubhead
Any $x\in V$ can be written uniquely in the form
$$x=e_{[a_1,b_1]}+e_{[a_2,b_2]}+\do+e_{[a_r,b_r]}$$
where $[a_r,b_r]\in\ci_D$ are such that any two of them are non-touching and
$r\ge0$, $1\le a_1\le b_1<a_1\le b_2<\do<a_r\le b_r\le D$.
Following \cite{\SYMP, 3.3} and \cite{\NEW, 1.11(a)} we set
$$u(v)=|\{s\in[1,r];a_s=_20,b_s=_21\}|-|\{s\in[1,r];a_s=_21,b_s=_20\}|\in\ZZ.\tag a$$
This defines a function $u:V@>>>\ZZ$. When $D\ge2$ we denote by $u':V'@>>>\ZZ$ the analogous
function with $D$ replaced by $D-2$. The following result appears also in 
\cite{\NEW, 1.11(b)}.

(b) {\it Assume that $D\ge2,i\in[1,D]$. Let $v'\in V'$ and let $v=T_i(v')+ce_i\in V$ where $c\in\FF_2$.
We have $u(v)=u'(v')$.}
\nl
We write $v'=e'_{[a'_1,b'_1]}+e'_{[a'_2,b'_2]}+\do+e'_{[a'_r,b'_r]}$ where $r\ge0$, 
$[a'_s,b'_s]\in\ci_{D-2}$ for all $s$ and any two of $[a'_s,b'_s]$ are non-touching. For each $s$, we 
have $T_i(e'_{[a'_s,b'_s]})=e_{[a_s,b_s]}$ where $[a_s,b_s]=\x_i[a'_s,b'_s]$ so that $a_s=_2a'_s$,
$b_s=_2b'_s$ and the various $[a_s,b_s]$ which appear are still non-touching with each other. Hence 
$u(T_i(v'))=u'(v')$. We have $v=T_i(v')$ or $v=T_i(v')+e_i$. If $v=T_i(v')$, we have 
$u(v)=u'(v')$, as desired. Assume now that $v=T_i(v')+e_i$.
From the definition of $\x_i$ we see that either

(i) $[i,i]$ is non-touching with any $[a_s,b_s]$, or

(ii) $[i,i]$ is not non-touching with some $[a,b]=[a_s,b_s]$ which is uniquely
determined and we have $a<i<b$.

If (i) holds then $e_i$ does not contribute to $u(v)$ and $u(v)=u(T_i(v'))=u'(v')$. 
We now assume that (ii) holds. Then $e_{[a,b]}+e_i=e_{[a,i-1]}+e_{[i+1,b]}$. We consider six cases.

(1) $a$ is even $b$ is odd, $i$ is even; then $|[i+1,b]|$ is odd so that the contribution of  
$e_{[a,i-1]}+e_{[i+1,b]}$ to $u(v)$ is $1+0$; this equals the contribution of $e_{[a,b]}$ to $u(T_i(v'))$ 
which is $1$.

(2) $a$ is even, $b$ is odd, $i$ is odd; then $|[a,i-1]|$ is odd so that the contribution of  
$e_{[a,i-1]}+e_{[i+1,b]}$ to $u(v)$ is $0+1$; this equals the contribution of $e_{[a,b]}$ to $u(T_i(v'))$ 
which is $1$.

(3) $a$ is odd, $b$ is even, $i$ is even; then $|[i+1,b]|$ is odd so that the contribution of  
$e_{[a,i-1]}+e_{[i+1,b]}$ to $u(v)$ is $0-1$; this equals the contribution of $e_{[a,b]}$ to $u(T_i(v'))$ 
which is $-1$.

(4) $a$ is odd, $b$ is even, $i$ is odd; then $|[a,i-1]|$ is odd so that the contribution of  
$e_{[a,i-1]}+e_{[i+1,b]}$ to $u(v)$ is $-1+0$; this equals the contribution of $e_{[a,b]}$ to $u(T_i(v'))$
 which is $-1$.

(5) $a=_2b=_2i+1$; then $|[a,i-1]|$ is odd, $|[i+1,b]|$ is odd so that the contribution of  
$e_{[a,i-1]}+e_{[i+1,b]}$ to $u(v)$ is $0+0$; this equals the contribution of $e_{[a,b]}$ to $u(T_i(v'))$ 
which is $0$.

(6) $a=_2b=_2i$; then the contribution of $e_{[a,i-1]}+e_{[i+1,b]}$ to $u(v)$ is $1-1$ or $-1+1$; this 
equals the contribution of $e_{[a,b]}$ to $u(T_i(v'))$ which is $0$.
\nl
This proves (b).

\mpb

Let $s\in[0,D/2]$ and let $x\in V(s)$, see 1.9. From the definition, the following holds:

(c) {\it If $s$ is even then $u(x)=s/2$; if $s$ is odd then $u(x)=-(s+1)/2$.}
\nl
We now define $\ti u:V@>>>\NN$ by $\ti u(x)=2u(x)$ if $u(x)\ge0$, $\ti u(x)=-2u(x)-1$ if $u(x)<0$. 
From (c) we deduce:

(d) {\it If $s\in[0,D/2]$ and $x\in V(s)$ then $\ti u(x)=s$.}

\subhead 1.12\endsubhead
As in \cite{\NEW, 1.12}. we
view $V$ as the set of vertices of a graph in which $x,x'$ in $V$ are joined whenever 
there exists $i\in[1,D]$ such that $x+x'=e_i$, $(x,e_i)=(x',e_i)=0$. (We then write $x\di x'$.) 
We show:

(a) {\it Let $s\in[0,D/2]$ and let $x,x'$ be in $V(s)$, see 1.9. Then $x,x'$ are in the same connected 
component of the graph $V$.}
\nl
As in 1.9 we set $t=s/2$ if $s$ is even and $t=(s+1)/2$ if $s$ is odd.
There is a unique element $x_s\in V(s)$ such that $n(x_s)\le n(y)$ 
for any $y\in V(s)$ (see 1.9 for the
the definition of $n(y)$). This element is of the form 
$x_s=e_{[a^0_1,b^0_1]}+e_{[a^0_2,b^0_2]}+\do+e_{[a^0_t,b^0_t]}$  where
$a^0_1,b^0_1,a^0_2,b^0_2,a^0_3,b^0_3,\do$ is $2,3,6,7,10,11,\do$ if $s$ is even and is
$1,2,5,6,9,10,\do$ if $s$ is odd. 
Let $\G$ be the connected component of the graph $V$ that contains $x_s$.
Let $x=e_{[a_1,b_1]}+e_{[a_2,b_2]}+\do+e_{[a_t,b_t]}\in V(s)$ be as in the definition of $V(s)$, see 1.9.
We show that $x\in\G$ by induction on $n(x)$. If $n(x)=n(x_s)$ then $x=x_s$ and there is nothing to prove.
Assume now that $n(x)>n(x_s)$. Then one of (i), (ii) below holds:

(i) for some $z\ge1$ we have $a_j=a_j^0$, $b_j=b_j^0$ for $j\in[1,z-1]$, $a_z>a_z^0$;

(ii) for some $z\ge1$ we have $a_j=a_j^0$, $b_j=b_j^0$ for $j\in[1,z-1]$, $a_z=a_z^0$, $b_z>b_z^0$.
\nl
In case (i) we have $a_z-2\in[1,2d]$, $(e_{a_z-2},x)=0$ hence $x+e_{a_z-2}\di x$. We have 
$(e_{a_z-1},x+e_{a_z-2})=0$ hence $x':=x+e_{a_z-2}+e_{a_z-1}\di x+e_{a_z-2}$. We have $n(x')=n(x)-2$. By 
the induction hypothesis we have $x'\in\G$ hence $x\in\G$.

In case (ii) we have $(e_{b_z-1},x)=0$ hence $x+e_{b_z-1}\di x$. We have $(e_{b_z},x+e_{b_z-1})=0$ hence 
$x':=x+e_{b_z-1}+e_{b_z}\di x+e_{b_z-1}$. We have $n(x')=n(x)-2$. By the induction hypothesis we have 
$x'\in\G$ hence $x\in\G$. This proves (a).

\subhead 1.13\endsubhead
For $x\in V$ we show:

(a) {\it there exists $s\in[0,D/2]$ and $\ti x\in V(s)$ (see 1.9) such that $x,\ti x$ are in the same 
component of the graph $V$.}
\nl
We argue by induction on $D$. If $D=0$ there is nothing to prove. Assume now that $D\ge1$. Assume first 
that $x$ is the element described in 1.9(d) or (e). Then $x\in V(D/2)$ so that there is nothing to prove.
Next we assume that $x$ is not the element described in 1.9(d) or (e). Then $(x,e_i)=0$ for some
$i\in[1,D]$. 
By 1.9(a) we have $x=T_i(x')+ce_i$ for some $x'\in V'$ and some $c\in\FF_2$. 
We first show the following result which appears also in \cite{\NEW, 1.12(a)}.

(b) {\it If $y,y'$ in $V'$ are joined in the graph $V'$ (analogue of the graph $V$)
 then $T_i(y),T_i(y')$ are in the same connected component of the graph $V$.}
\nl
We can find $j\in[1,D-2]$ such that $(y,e'_j)'=(y',e'_j)'=0$, $y+y'=e'_j$.
Hence $(\ti y,T_i(e'_j))=(\ti y',T_i(e'_j))=0$, $\ti y+\ti y'=T_i(e'_j)$
where $\ti y=T_i(y),\ti y'=T_i(y')$.
If $T_i(e'_j)=e_h$ for some $h\in[1,D]$ then $\ti y,\ti y'$ ar joined in $V$, as required.
If this condition is not satisfied then $1<i<D$, $j=i-1$ and $T_i(e'_j)=e_j+e_{j+1}+e_{j+2}$.
We have $(\ti y,e_j+e_{j+1}+e_{j+2})=0$,  $\ti y+\ti y'=e_j+e_{j+1}+e_{j+2}$. Since 
$\ti y\in V_i$ we have $(\ti y,e_i)=0$ hence $(\ti y,e_{j+1})=0$ so that
$(\ti y,e_j)=(\ti y,e_{j+2})$. We are in one of the two cases below.

(1) We have $(\ti y,e_j)=(\ti y,e_{j+2})=0$.

(2) We have $(\ti y,e_j)=(\ti y,e_{j+2})=1$.
\nl
In case (1) we consider the four term sequence 
$\ti y,\ti y+e_j,\ti y+e_j+e_{j+2},\ti y+e_j+e_{j+1}+e_{j+2}=\ti y'$; any two
consecutive terms of this sequence are joined in the graph $V$.
In case (2) we consider the four term sequence $\ti y,\ti y+e_{j+1},\ti y+e_j+e_{j+1},
\ti y+e_j+e_{j+1}+e_{j+2}=\ti y'$; any two
consecutive terms of this sequence are joined in the graph $V$. We see that in both cases $\ti y,\ti y'$ 
are in the same connected component of $V$ and (b) is proved.

We now continue the proof of (a). By the induction hypothesis there exists $s\in[0,(D/2)-1]$ and
$x''\in V'(s)$ such that $x',x''$ are in the same connected component of $V'$.
Here $V'(s)$ is defined like $V(s)$ (replacing $V$ by $V'$).
By (b), $T_i(x'),T_i(x'')$ are in the same connected component of $V$. 
From the definitions we see that $T_i(V'(s))\sub V(s)$. Thus $T_i(x'')\in V(s)$.
Clearly $x,T_i(x')$ are joined in the graph $V$. Hence $x,T_i(x'')$ are joined in the graph $V$. 
We see that (a) holds.

\subhead 1.14\endsubhead
The following result follows by repeated application of 1.11(b).

(a) {\it If $x,x'$ in $V$ are in the same connected component of the graph $V$ then $u(x)=u(x')$.}
\nl
We can assume that $x,x'$ are joined in the graph $V$. Then for some $i\in[1,D]$ 
we have $x=T_i(y)+ce_i$, $x'=T_i(y)+c'e_i$ where $y\in V'$, $c\in\FF_2,c'\in\FF_2$.
By 1.11(b) we have $u(x)=u'(y)$, $u(x')=u'(y)$, hence  $u(x)=u(x')$. This proves (a).

We now show the converse.

(b) {\it If $x,x'$ in $V$ satisfy $u(x)=u(x')$, then $x,x'$ are in the same connected component of the graph
 $V$.}
\nl
By 1.13(a) we can find $s,s'$ in $[0,D/2]$ and $x_1\in V(s)$, $x'_1\in V(s')$ such that $x,x_1$ are in the 
same connected component of the graph $V$ and $x',x'_1$ are in the same connected component of the graph 
$V$. Thus, it is enough to prove that $x_1,x'_1$ are in the same connected component of the graph $V$. By 
(a), we have $u(x_1)=u(x'_1)$ hence $\ti u(x_1)=\ti u(x'_1)$. From 1.11(d) we have 
$\ti u(x_1)=s$, $\ti u(x'_1)=s'$. Using $\ti u(x_1)=\ti u(x'_1)$ we deduce that $s=s'$. Since 
$x_1\in V(s),x'_1\in V(s)$, they are in the same connected component of the graph $V$, by 1.12(a). This 
proves (b).

\mpb

We show:

(c) {\it Let $B\in\SS_D$. Let $k=|B^0|$, $k'=\ti u(\e(B))\in\ZZ$. Then $k'=k$.}
\nl
We argue by induction on $D$. If $D=0$ there is nothing to prove. Assume now that $D\ge2$.
If $B\in\SS_D^{prim}$, then $\e(B)\in V(k)$, see 1.9(g), and the result follows from 1.11(d).
We now assume that $B\n\SS_D^{prim}$. We can find $i\in[1,D]$ and $B'\in\SS_{D-2}$ such that
$B=t_i(B')$. By 1.9(b) we have $\e(B)=T_i(\e'(B'))+ce_i$ where $c\in\FF_2$.
Using 1.11(b) we deduce $u(\e(B))=u'(\e'(B'))$. Hence $\ti u(\e(B))=\ti u'(\e'(BB'))$ where 
$\ti u':V'@>>>\ZZ$ is defined in terms of $u'$ in the same way as $\ti u$ is defined in terms of $u$. By
the induction hypothesis we have $\ti u'(\e'(B'))=k$. This proves (c).

\mpb

(d) {\it The map $\e:\SS_D@>>>V$, see 1.9, is injective.}
\nl
Assume that $B,\tB$ in $\SS_D$ are such that $\e(B)=\e(\tB)$. Let $k=|B^0|$, $\ti k=|\tB^0|$. Let 
$k'=\ti u(\e(B)))=\ti u(\e(\tB))$. By (c) we have $k'=k$, $k'=\ti k$. It follows that $k=\ti k$.
Using now 1.10(a), we see that $B=\tB$. This proves (d).

\subhead 1.15\endsubhead
Let $s\in[0,D/2]$ and let $B\in\SS_D$ be as in 1.9(f), so that $\e(B)\in V(s)$. We show:

(a) {\it For any $x\in\la B\ra$ we have $\ti u(x)\le s$; moreover, we have $\ti u(x)=s$ for a unique 
$x\in\la B\ra$.}
\nl
We argue by induction on $D$. If $D=0$ the result is obvious. We now assume that $D\ge2$.
Assume first that $s=D/2$. If $x_1=\e(B)$ then $x_1\in\la B\ra$, see 1.9(c), and $\ti u(x_1)=D/2$, see 1.14(c).
Conversely, assume that $x'\in V$, $\ti u(x')=D/2$. Using 1.14(b), we see that $x',x_1$ are in the same 
connected component of $V$. From 1.9(d),(e), we see that $(x_1,e_i)=1$ for any $i\in[1,D]$. Thus,
$x_1$ is a connected component of $V$ by itself, so that $x'=x_1$. Hence in this case (a) holds.
Next we assume that $s<D/2$. Then $B'=\{[1,D-2],[2,D-3],\do,[s,D-1-s]\}\in\SS_{D-2}$ satisfies 
$\e'(B')=s$ (by 1.9(f)). We have $B=\{\x_i(I');I'\in B'\}$. Let $i=s+1$. Let $\tB=t_i(B')=B\sqc\{i\}$.
Using the induction hypothesis for $B'$ and 1.11(b) we see that for any 
$$x\in T_i(\la B'\ra)\op\FF_2e_i=\la B\ra\op\FF_2e_i=\la\tB\ra$$ 
(see 1.9(c)) we have $\ti u(x)\le s$; 
moreover, we have $\ti u(x)=s$ for exactly two values of $x\in\la\tB\ra$ (whose sum is $e_i$). One of 
these values is in $\la B\ra$ and the other is not in $\la B\ra$. This proves (a).

\subhead 1.16\endsubhead
Let $F$ be the $\CC$-vector space consisting of functions $V@>>>\CC$. For $x\in V$ let $\ps_x\in F$ be the
characteristic function of $x$. For $B\in\SS_D$ let $\Ps_B\in F$ be the characteristic function of 
$\la B\ra$. Let $\tF$ be the $\CC$-subspace of $F$ generated by $\{\Ps_B;B\in\SS_D\}$. 
When $D\ge2$ we define $\ps'_{x'}$ for $x'\in V'$ and $\Ps'_{B'}$ for $B'\in\SS_{D-2}$, $F',\tF'$,
in terms of $\SS_{D-2}$ in the same way as $\ps_x,\Ps_B,F,\tF$ were defined in terms of $\SS_D$. 
For any $i\in[1,D]$ we define a linear map $\th_i:F'@>>>F$ by $f'\m f$ where 
$f(T_i(x')+ce_i)=f'(x')$ for $x'\in V',c\in\FF_2$, $f(x)=0$ for $x\in V-e_i^\pe$. We have 

$\th_i(\ps'_{x'})=\ps_{T_i(x')}+\ps_{T_i(x')+e_i}$ for any $x'\in V'$,

$\th_i(\Ps'_{B'})=\Ps_{t_i(B')}$ for any $B'\in\SS_{D-2}$. 
\nl
We show:

(a) {\it For any $x\in V$, we have $\ps_x\in\tF$.}
\nl
We argue by induction on $D$. If $D=0$ the result is obvious. We now assume that $D\ge2$. We first show:

(b) {\it If $x,\tx$ in $V$ are joined in the graph $V$ and if (a) holds for $x$, then (a) holds for $\tx$.}
\nl
We can find $j\in[1,D]$ such that $x+\tx=e_j$, $(x,e_j)=0$.
We have $x=T_j(x')+ce_j$, $\tx=T_j(x')+c'e_j$ where $x'\in V'$ and $c\in\FF_2,c'\in\FF_2$, $c+c'=1$.
 By the induction 
hypothesis we have $\ps'_{x'}=\sum_{B'\in\SS_{D-2}}a_{B'}\Ps'_{B'}$ where $a_{B'}\in\CC$. 
Applying $\th_j$ we obtain
$$\ps_x+\ps_{\tx}=\sum_{B'\in\SS_{D-2}}a_{B'}\Ps_{t_j(B')}$$ 
We see that $\ps_x+\ps_{\tx}\in\tF$. Since $\ps_x\in\tF$, by assumption, we see that $\ps_{\tx}\in\tF$. 
This proves (b).

\mpb

For any $s\in[0,D/2]$ we show:

(c) {\it If $x\in V$ is such that $\ti u(x)=s$ then $\ps_x\in\tF$.}
\nl
We argue by induction on $s$. Let $\tx=\e(B)$ where $B$ is as in 1.9(f) so that $\tx\in V(s)$ (see 1.9(g))
and $\ti u(\tx)=s$ (see 1.11(d)). Using 1.14(b) we see that $x,\tx$ are in the same connected
component of the graph $V$ and using (b) see that it is enough to show that $\ps_{\tx}\in\tF$.
Let $x_0$ be the unique element of $\la B\ra$ such that $\ti u(x_0)=s$ (see 1.15(a)). By the uniqueness 
of $x_0$ we must have $x_0=\tx$. From 1.15 we see that for $x_1\in\la B\ra-\{\tx\}$ we have 
$\ti u(x_1)<s$; for such $x_1$ we have $\ps_{x_1}\in\tF$ by the induction hypothesis. We have 
$\Ps_B=\ps_{\tx}+\sum_{x_1\in\la B\ra-\{\tx\}}\ps_{x_1}=\ps_{\tx}\mod\tF$. Since $\Ps_B\in\tF$, we see that
$\ps_{\tx}\in\tF$. This proves (c) hence also (a).

\mpb

Since $\tF\sub F$, we see that (a) implies:

(d) $F=\tF$.
\nl
This extends \cite{\NEW, 1.15(c)}. We have the following result which extends
\cite{\NEW, 1.16}.

\proclaim{Theorem 1.17} (a) $\{\Ps_B;B\in\SS_D\}$ is a $\CC$-basis of $F$.

(b) $\e:\SS_D@>>>V$ is a bijection.
\endproclaim
From the definition of $\tF$ we have $\dim\tF\le|\SS_D|$. By 1.14(d) we have $|\SS_D|\le|V|=\dim F$.
Since $F=\tF$ (see 1.16(d)), it follows that $\dim\tF=|\SS_D|=|V|=\dim F$. Using again the definition
of $\tF$ and the equality $F=\tF$ we see that (a) holds. Since the map in (b) is injective (see 1.14(d))
and $|\SS_D|=|V|$ we see that it is a bijection so that (b) holds.

\mpb

Let $\pmb\cf(V)$ be the set of (isotropic) subspaces of $V$ of the form $\la B\ra$ for some $B\in\SS_D$.
By definition, the map $\SS_D@>>>\pmb\cf(V)$, $B\m\la B\ra$ is surjective. In fact,

(c) {\it this map is a bijection.}
\nl
Indeed, if $B,\tB$ in $\SS_D$ satisfy $\la B\ra=\la\tB\ra$ then the functions $\Ps_B,\Ps_{\tB}$ in $F$
coincide and (d) follows from (a).

\mpb

Note that $\pmb\cf(V)$ admits an inductive definition similar to that of $\SS_D$. If $D=0$, $\pmb\cf(V)$ consists of
the subspace $\{0\}$. If $D\ge2$, a subspace $E$ of $V$
is in $\pmb\cf(V)$ if and ony if it is either of the form $\la B\ra$ for some $B\in\SS_D^{prim}$ or if
there exists $i\in[1,D]$ and $E'\in\pmb\cf(V')$ such that $E=T_i(E')\op\FF_2 e_i$. 

\subhead 1.18\endsubhead
Assume that $D\ge2$. Let $B\in\SS_D$ and let $i\in[1,D]$ be such that $\{i\}\in B$.  
Let $Z_i$ be the set of all $[a,b]\in B^1$ such that $a<i<b$. If $I\in Z_i,I'=[a',b']\in Z_i$, then
$I\cap I'\ne\emp$ hence we have either $I\sub I'$ or $I'\sub I$. It follows that if $Z_i\ne\emp$ then 
$Z_i$ contains a unique interval $[a,b]$ such that $b-a$ is minimum; we set $Z_i^{min}=\{[a,b]\}$. We show:

(a) {\it If $Z_i\ne\emp$ and $Z_i^{min}=\{[a,b]\}$, then $a=_2b=_2i+1$.} 
\nl
If this is not so, then $a=_2b=_2i$. By $(P_1)$ there exists $[a_1,b_1]\in B^1$ such that 
$a<a_1\le i-1\le b_1<b$. If $b_1=i-1$, then applying $(P_0)$ to $[a_1,b_1],\{i\}$ gives a contradiction. 
Thus $b_1\ge i$ and $i\in[a_1,b_1]$. By the minimality of $b-a$, we have $[a_1,b_1]=\{i\}$. This 
contradicts $i-1\in[a_1,b_1]$ and proves (a).

\mpb

Let $h_0<h_1<\do<h_{2k+1}$ be the sequence attached to $B$ in $(P_2)$. We show:

(b) {\it Assume that $h_s<i<h_{s+1}$. If $s\in[0,k-1]$ and $i=_2s$, then $Z_i\ne\emp$.
If $s\in[k+1,2k]$ and $i=_2s+1$, then $Z_i\ne\emp$.}
\nl
We prove the first assertion of (b).
We have $h_s<i-1<h_{s+1}$ (since $h_s\ne_2 i-1$). By $(P_2)$ we can find $[a,b]\in B^1$ such that 
$h_s<a\le i-1\le b<h_{s+1}$. If $b=i-1$ then applying $(P_0)$ to $[a,b],\{i\}$ gives a contradiction. 
Thus, $b\ge i$ and $i\in[a,b]$. Since $a<i$ we have $\{i\}\prec[a,b]$ so that $[a,b]\in Z_i$. This 
proves the first assertion of (b). The second assertion of (b) can be deduced from the first assertion
using the involution $\t_D:\SS_D@>>>\SS_D$ in 1.2.

\mpb

We show:

(c) {\it If $h_s<i<h_{s+1}$, $s\in[0,k-1]$, $i=_2s+1$, then either $Z_i=\emp$ or $Z_i\ne\emp$ and
$Z_i-Z_i^{min}\ne\emp$.}
\nl
Assume that $Z_i\ne\emp$. Let $[a,b]\in Z_i^{min}$, so that $a<i<b$. Using 1.3(e) we see that 
$h_s<a<b<h_{s+1}$. By (a) we have $a=_2h_s$. Since $h_s<a$, we must have $h_s<a-1<h_{s+1}$. By 
$(P_2)$ we can find $[a',b']\in B^1$ such that $h_s<a'\le a-1\le b'<h_{s+1}$. If $b'=a-1$ then 
applying $(P_0)$ to $[a,b],[a',b']$ gives a contradiction. Thus, $b'\ge a$, so that 
$[a',b']\cap[a,b]\ne\emp$. This implies that either $[a',b']\sub[a,b]$ or $[a,b]\prec[a',b']$. The first 
alternative does not hold since $a-1\in[a',b'],a-1\n[a,b]$. Thus we have $[a,b]\prec[a',b']$ so that 
$[a',b']\in Z_i-Z_i^{min}$. This proves (c).

\mpb

We define a collection $C$ of subsets of $\ci_D$ as follows:

(i) If $h_k<i<h_{k+1}$ and $Z_i=\emp$ then $C=B-\{i\}$.

(ii) If $h_k<i<h_{k+1}$ and $Z_i\ne\emp$ then $C=(B-\{[a,b],\{i\}\})\sqc\{[a,i-1],[i+1,b]\}$
where $Z_i^{min}=\{[a,b]\}$. 

(iii) If $h_s<i<h_{s+1}$, $s\in[0,k-1]$, $i=_2s$, so that $Z_i\ne\emp$ (see (b)) then
$C=(B-\{[a,b],\{i\}\})\sqc\{[a,i-1],[i+1,b]\}$ where $Z_i^{min}=\{[a,b]\}$. 

(iv) If $h_s<i<h_{s+1}$, $s\in[0,k-1]$, $i=_2s+1$ and $Z_i\ne\emp$ then 
$C=(B-\{[a,b],\{i\}\})\sqc\{[a,i-1],[i+1,b]\}$ where $Z_i^{min}=\{[a,b]\}$. 

(v)  If $h_s<i<h_{s+1}$, $s\in[k+1,2k]$, $i=_2s+1$ so that $Z_i\ne\emp$ (see (b)) then
$C=(B-\{[a,b],\{i\}\})\sqc\{[a,i-1],[i+1,b]\}$ where $Z_i^{min}=\{[a,b]\}$. 

(vi) If $h_s<i<h_{s+1}$, $s\in[k+1,2k]$, $i=_2s$ and $Z_i\ne\emp$ then 
$C=(B-\{[a,b],\{i\}\})\sqc\{[a,i-1],[i+1,b]\}$ where $Z_i^{min}=\{[a,b]\}$. 

(vii) If $h_s<i<h_{s+1}$, $s\in[0,k-1]$, $i=_2s+1$ and $Z_i=\emp$ then 
$C=(B-\{[h_{s+1},h_{2k-s}],\{i\}\})\sqc\{[i,h_{2k-s}],[i+1,h_{s+1}-1]\}$.

(viii) If $h_s<i<h_{s+1}$, $s\in[k+1,2k]$, $i=_2s$ and $Z_i=\emp$ then 
$C=(B-\{[h_{2k-s+1},h_s],\{i\}\})\sqc\{[h_{2k-s},i],[h_s+1,i-1]\}$.

\mpb

For $h\in\{0,1\}$ let $C^h$ be the set of all $[a',b']\in C$ such that $b-a=_2h+1$. We show:

(d) {\it $C$ satisfies properties $(P_0),(P_1),(P_2)$.}
\nl
We refer to properties $(P_0),(P_1),(P_2)$ for $C$ as $(P'_0),(P'_1),(P'_2)$.
The verification of $(P'_0)$ is immediate. We check $(P'_2)$. The sequence $h'_0<h'_1<\do<h'_{2k+1}$ in 
$(P'_2)$ is:

$h_0<h_1<\do<h_{2k+1}$ (of $(P_2)$ for $B$) in cases (i)-(vi) (in these cases we use that 
$a=_2i+1,b=_2i+1$, see (a)); 

$h_0<h_1<\do<h_s<i<h_{s+2}<\do<h_{2k+1}$ in case (vii);

$h_0<h_1<\do<h_{s-1}<i<h_{s+1}<\do<h_{2k+1}$ in case (viii).

We check $(P'_1)$. In case (i), $(P'_1)$ is immediate.
In case (ii)-(vi) let $c$ be such that $a<c<i-1$ or $i+1<c<b$, $c=_2a+1$. By $(P_1)$ for $B$ we can find 
$[a_1,b_1]\in B^1$ such that $a<a_1\le c\le b_1<b$. If $c<i-1$, $b_1\ge i$ then $[a_1,b_1]\in Z_i$, 
contradicting $Z_i=\emp$; if $i+1<c$, $a_1\le i$, then $[a_1,b_1]\in Z_i$, contradicting $Z_i=\emp$. Thus,
we have $a<a_1\le c\le b_1\le i-1$ or $i+1\le a_1\le c\le b_1<b$. If $b_1=i-1$ or $a_1=i+1$,  
then applying $(P_0)$ for $B$ to $[a_1,b_1],\{i\}$ gives a contradiction;
thus we have $a<a_1\le c\le b_1<i-1$ or $i+1<a_1\le c\le b_1<b$. Moreover, since $[a_1,b_1]\in B^1$ we 
have $[a_1,b_1]\in C^1$ so that $(P'_1)$ holds.
In case (vii) let $c$ be such that $i+1<c<h_{s+1}-1$, $c=_2i$. By $(P_2)$ for $B$ we can find 
$[a,b]\in B^1$ such that $h_s<a\le c\le b<h_{s+1}$. We have $b\le h_{s+1}-1$. If $a\le i$, then 
$[a,b]\in Z_i$, contradicting $Z_i=\emp$. Thus, $a>i$, so that $i+1\le a\le c\le b\le h_{s+1}-1$. If 
$b=h_{s+1}-1$ then applying $(P_0)$ for $B$ to $[a,b],[h_{s+1},h_{2k-s}]$ gives a contradiction.
Thus, $b<h_{s+1}-1$. If $a=i+1$ then applying $(P_0)$ for $B$ to $[a,b],\{i\}$ gives a contradiction.
Thus $i+1<a$. Moreover, since $[a,b]\in B^1$ we have $[a,b]\in C^1$ so that $(P'_1)$ holds.
In case (viii), $(P'_1)$ is proved by an argument similar (and symmetric under $\t_D$) to that in case 
(vii).

We check $(P'_2)$ with $j\in[0,k-1]$. In case (i), $(P'_2)$ is immediate. Let $c$ be such that 
$h'_j<c<h'_{j+1}$, $c=_2j+1$.  In cases (ii)-(vi), by $(P_2)$ for $B$ we can find
$[a',b']\in B^1$ such that $h_j<a'\le c\le b'<h_{j+1}$. If we are in case (ii),(v) or (vi),
or (iii),(iv) with $s\ne j$, we have $[a',b']\in C^1$ and $(P'_2)$ holds. Assume that we are in case
(iii) or (iv) with $s=j$. Let $[a,b]\in B^1$ be the unique interval in $Z_i^{min}$.
If $[a',b']\ne[a,b]$, then $[a',b']\in C^1$ and $(P'_2)$ holds. Thus we can assume that $[a',b']=[a,b]$ so
that $a\le c\le b$. If $i\n[a',b']$ then $[a',b']\in C^1$ and $(P'_2)$ holds. Thus we can assume that 
$i\in[a,b]=[a',b']$. In case (iii) (with $s=j$) we have $c\ne i$ (since $c=_2j+1,i=_2s,s=_2j$) 
hence $c<i$ or $c>i$. Thus we have $c\in[a,i-1]$ or $c\in[i+1,b]$ and $[a,i-1]\in C^1$, $[i+1,b]\in C^1$ 
and $(P'_2)$ holds. In case (iv) with $s=j$, by (c) we can find $[a'',b'']\in Z_i$ such that 
$[a,b]\prec[a'',b'']$. We have $[a'',b'']\in C^1$ and $h_j<a''\le c\le b''<h_{j+1}$. Thus, $(P'_2)$ holds.
Assume now that we are in case (vii). If $j\ne s+1$ then by $(P_2)$ for $B$ we can find 
$[a,b]\in B^1$ such that $h_j<a\le c\le b<h_{j+1}$.
If in addition we have $j\ne s$ then $h'_j<a\le c\le b<h'_{j+1}$, $[a,b]\in C^1$ and $(P'_2)$ holds. 
If $j=s$, we jave $c<i$ hence $a<i$. We show that $h_s<a\le c\le b<i$ (in particular, $[a,b]\in C^1$). 
Now $h_s<a$ holds since $h_s=h'_j$. To prove that $b<i$, we assume that $i\le b$ so that 
$i\in[a,b]$. Since $Z_i=\emp$ we deduce that $a=b=i$ hence $c=i$. This contradicts
$c<h'_{s+1}=i$ and proves $(P'_2)$ in this case. 
If $j=s+1$, then taking $[a,b]=[i+1,h_{s+1}-1]\in C^1$, we have 
$h'_{s+1}<i+1\le c\le h_{s+1}-1<h'_{s+2}$ so that $(P'_2)$ holds.  

Assume now that we are in case (viii). By $(P_2)$ for $B$ we can find 
$[a,b]\in B^1$ such that $h_j<a\le c\le b<h_{j+1}$ hence $h'_j<a\le c\le b<h'_{j+1}$. We have 
$[a,b]\in C^1$ so that $(P'_2)$ holds.

The proof of $(P'_2)$ with $j\in[k+1,2k]$ is similar (and symmetric under $\t_D$) to the  
proof of $(P'_2)$ with $j\in[0,k-1]$. This completes the proof of (d).

\mpb

From (d) and 1.3(c) we deduce:

(e) {\it We have $C\in\SS_D$.}
\nl
From the definitions we deduce:

(f) {\it For $j\in[1,D]-\{i\}$ we have $f_j(C)=f_j(B)$. In case (i) we have $f_i(C)=f_i(B)-1$.
In cases (ii)-(viii) we have $f_i(C)=f_i(B)-2$.}
\nl
From (f) we deduce:

(g) {\it For $j\in[1,D]-\{i\}$ we have $\e_j(C)=\e_j(B)$. We have $\e_i(C)=\e_i(B)+1$.}
\nl
(For the second assertion of (g) in case (ii)-(viii) we use 1.6(d); in case (i) we have 
$f_i(C)=-k-\un k,f_i(B)=-k-\un k+1$ and $k+\un k\in2\ZZ$, so that the second assertion of (g) holds by 
1.6(e).)

We show:

(h) {\it We have $\e(C)=\e(B)+e_i$, $\e(B)\in e_i^\pe$. In other words, $\e(C),\e(B)$ are
joined in the graph $V$.}
\nl
The first assertion of (h) is a restatement of (g). For the second assertion we note that by 1.3(f)
we have $B=t_i(B')$ for some $B'\in\SS_{D-2}$, so that $\la B\ra\sub V_i\op\FF_2e_i$ and it remains to use
1.9(a) and 1.9(c).

\mpb

(i) {\it We shall also use the notation $C=B[i]$ when $C$ is os obtained from $B,i$ as above.}

\subhead 1.19\endsubhead
We view $\SS_D$ as the set of vertices of a graph in which $B_1,B_2$ in $\SS_D$ are joined whenever
$\e(B_1)\di\e(B_2)$, see 1.12. (We then write $B_1\di B_2$.) Thus the bijection $\e:\SS_D@>\si>>V$ is a graph 
isomorphism. We show:

(a) {\it Let $B_1,B_2$ in $\SS_D$ be such that $B_1\di B_2$. Define $i\in[1,D]$ by $e_i=\e(B_1)+\e(B_2)$.
Then $\{i\}$ belongs to exactly one of $B_1,B_2$, say $B_1$ and we have $B_2=B_1[i]$, see 1.18(i).
Moreover, we have $B_1=t_i(B')$ for a well defined $B'\in\SS_{D-2}$ and $\e(B_1)=T_i(\e'(B'))\mod\FF_2e_i$, 
$\e(B_2)=T_i(\e'(B'))\mod\FF_2e_i$.}
\nl
We have $\e(B_1)=T_i(x')+c_1e_i$, $\e(B_2)=T_i(x')+c_2e_2$ for a well defined $x'\in V'$, $c_1\in\FF_2$,
$c_2\in\FF_2$ such that $c_1+c_2=1$.
Define $B'\in\SS_{D-2}$ by $\e'(B')=x'$. By 1.9(b) we have $\e(t_i(B'))=T_i(x')+ce_i$ with $c\in\FF_2$.
Since $c_1+c_2=1$ we have $c=c_1$ or $c=c_2$. Assume for example that $c=c_1$. Then $\e(t_i(B'))=\e(B_1)$.
Since $\e$ is a bijection we deduce that $B_1=t_i(B')$, so that $\{i\}\in B_1$. Let $C_1=B_1[i]\in\SS_D$,
see 1.18(i). By 1.18(h) we have 
$\e(C_1)=\e(B_1)+e_i$ so that $\e(C_1)=T_i(x')+c_1e_i+e_i=T_i(x')+c_2e_2=\e(B_2)$. Since $\e$ is a 
bijection we deduce that $C_1=B_2$. Note that $\{i\}\n C_1$ so that $\{i\}\n B_2$. This proves (a).

\subhead 1.20\endsubhead
For $B,\tB$ in $\SS_D$ we say that $B\le\tB$ if either 

(i) $|B^0|<|\tB^0|$ or 

(ii) $|B^0|=|\tB^0|$ and for any $i\in[1,D]$ we have $f_i(B)\le f_i(\tB)|$. 
\nl
We show:

(a) {\it This is a partial order on $\SS_D$.}
\nl
It is enough to prove that for $B,\tB$ in $\SS_D$ such that $B\le\tB$ and $\tB\le B$ we have $B=B'$. We have 
$|B^0|\le|\tB^0|\le|B^0|$ hence $|B^0|=|\tB^0|$ and $f_i(B)=f_i(\tB)$ for all $i$ hence 
$\e_i(B)=\e_i(\tB)$ and $\e(B)=\e(\tB)$. Since $\e$ is a bijection (1.17(b)), we deduce that $B=B'$. 
This proves (a).

\mpb

For $x,\tx$ in $V$ we say that $x\le\tx$ if $\e\i(x)\le\e\i(\tx)$ where $\e\i:V@>>>\SS_D$ is the
bijection inverse to $\e:\SS_D@>>>V$. This is a partial order on $V$. 
We shall write $x<\tx$ whenever $x\le\tx$ and $x\ne\tx$. Using the definitions and 1.4(a) we deduce:

(b) {\it Assume that $D\ge2$, $i\in[1,D]$, $B'\in\SS_{D-2}$, $\tB'\in\SS_{D-2}$. If $B'\le\tB'$, then 
$t_i(B)\le t_i(\tB)$. Hence if $x'\in V'$, $\tx'\in V'$, $x'\le\tx'$, then
$t_i(\e'{}\i(x')\le t_i(\e'{}\i(\tx')$.}

\mpb

Clearly, for any $x\in V$ we have $0\le x$. We denote by $\nu(x)$ the largest number $r\ge0$ such that 
there exists a sequence $0=x_0<x_1<\do<x_r=x$ in $V$. We have $\nu(0)=0$ and $\nu(x)>0$ if $x\ne0$.

\mpb

We show:

(c) {\it Assume that $B\in\SS_D^{prim}$. Recall that $z:=\e(B)\in\la B\ra$ (see 1.9(c)). If 
$y\in\la B\ra$ and $y\ne z$ then $y<z$.}
\nl
We set $k=|B^0|\in[0,D/2]$. By 1.9(g) we have $z\in V(k)$ and by 1.15(a) we have $\ti u(x_0)=k$ for a 
unique $x_0\in\la B\ra$, $\ti u(x)<k$ for any $x\in\la B\ra$ such that $x\ne x_0$. By 1.14(c) we have 
$\ti u(z)=|B^0|=k$ so that $x_0=z$. Thus for $y$ as in (c) we have $y\ne x_0$ so that 
$\ti u(y)<k$ that is $\ti u(\e(B'))<k$ where $B'=\e\i(y)$. By 1.14(c) this implies $|B'{}^0|<k$ that is 
$|B'{}^0|<|B^0|$ so that $B'<B$ and $y<z$. This proves (c).

\mpb

Let $x\in V$. By 1.16(d) we have $\ps_x=\sum_{\tx\in V}c_{x,\tx}\Ps_{\e\i(\tx)}$ where 
$c_{x,\tx}\in\CC$. Moreover, by 1.17, the coefficients $c_{x,\tx}$ are uniquely determined. We 
state:

\proclaim{Theorem 1.21} If $x\in V,\tx\in V$, $c_{x,\tx}\ne0$, then $\tx\le x$. Moreover, $c_{x,x}=1$.
\endproclaim
We argue by induction on $D$; for fixed $D$ we argue by (a second) induction on $\nu(x)$. If $D=0$ the 
result is obvious. Now assume that $D\ge2$. Assume first that $\e\i(x)\in\SS_D^{prim}$. Since 
$x\in\la\e\i(x)\ra$, we have $\Ps_{\e\i(x)}=\ps_x+\sum_{x_1\in Z}\ps_{x_1}$ where 
$Z=\la\e\i(x)\ra-\{x\}$. By 1.20(c), for any $x_1\in Z$ we have $x_1<x$ so that $\nu(x_1)<\nu(x)$. 
By the (second) induction hypothesis, for any $x_1\in Z$, $\ps_{x_1}$ is a linear combination of 
elements $\Ps_{\e\i(x_2)}$ with $x_2\in V$, $x_2\le x_1$ (hence $x_2<x$). It follows that the statement 
of the theorem holds for our $x$.

Next we assume that $B=\e\i(x)\n\SS_D^{prim}$. We can find $i\in[1,D]$ such that $\{i\}\in\e\i(x)$.
We have $\e\i(x)=t_i(B')$ where $B'\in\SS_{D-2}$. Let $x'=\e'(B')\in V'$. We have $t_i(\e'{}\i(x'))=B$.
From the first induction hypothesis we have
$$\ps'_{x'}=\sum_{\tx'\in V';\tx'\le x'}c'_{x',\tx'}\Ps'_{\e'{}\i(\tx')}\tag a$$ 
where $c'_{x',\tx'}\in\CC$ and $c'_{x',x'}=1$.
Let $C=B[i]$, see 1.18(i). We have $|C^0|=|B^0|$ and from 1.18(f) we see that $C<B$ hence $y<x$
where $y=\e(C)\in V$. Applying to (a) $\th_i$ (as in the proof of 1.16(b)) we obtain
$$\ps_x+\ps_y=\sum_{\tx'\in V';\tx'\le x'}c'_{x',\tx'}\Ps_{t_i(\e'{}\i(\tx'))}$$ 
(we have used 1.19(a)). By 1.20(b) the inequality 
$\tx'\le x'$ implies $t_i(\e'{}\i(\tx'))\le t_i(\e'{}\i(x'))=B$;
moreover if  $\tx'\ne x'$ then $t_i(\e'{}\i(\tx'))\ne t_i(\e'{}\i(x'))=B$.
We see that $\ps_x+\ps_y$ is a linear combination of terms 
$\Ps_{\e\i(z)}$ with $z\in V$, $z\le x$, and the coefficient of $\Ps_{\e\i(x)}$ is $1$.

Since $y<x$ we have $\nu(y)<\nu(x)$. By the (second) induction hypothesis
$\ps_y$ is a linear commbination of terms $\Ps_{\e\i(z)}$ with $z\in V$, $z\le y$
hence $z<x$. We see that $\ps_x$ is a linear combination 
of terms $\Ps_{\e\i(z)}$ with $z\in V$, $z\le x$ and the coefficient of $\Ps_{\e\i(x)}$ is $1$.
This proves the theorem.

\subhead 1.22\endsubhead
For $x\in V$ we have $\Ps_{\e\i(x)}=\sum_{\tx\in V}d_{x,\tx}\ps_{\tx}$ where $d_{x,\tx}=1$ if 
$\tx\in\la\e\i(x)\ra$ and $d_{x,\tx}=0$ if $\tx\n\la\e\i(x)\ra$. Recall that $d_{x,x}=1$. We show:

(a) {\it If $d_{x,\tx}\ne0$ then $\tx\le x$.}
\nl
From the definitions for $x,x'$ in $V$ we have $\sum_{\tx\in V}c_{x,\tx}d_{\tx,x'}=\d_{x,x'}$ (Kronecker $\d$). 
Using 1.21 we deduce $d_{x,x'}+\sum_{\tx\in S_D;\tx<x}c_{x,\tx}d_{\tx,x'}=\d_{x,x'}$.
From this the desired result follows by induction on $\nu(x)$.

\mpb

We show:

(b) {\it There is a unique bijection $e:V@>\si>>\pmb\cf(V)$ (see 1.17) such that $x\in e(x)$ for any $x\in V$.}
\nl
The map $e:x\m\la \e\i(x)\ra$, $V@>>>\pmb\cf(V)$ is a well defined bijection, see 1.17(b),(c). For $x\in V$
we have $x\in e(x)$ by 1.9(c). This proves the existence of $e$. We prove uniqueness. Let
$e':V@>>>\pmb\cf(V)$ be a bijection such that $x\in e'(x)$ for any $x\in V$. 
We define a bijection $\s:V@>\si>> V$ by $\s=e'{}\i e$. Then for any $X\in\pmb\cf(V)$ we have
     $\s(e\i(X))=e'{}\i (X)$. Setting $x=e\i(X)$ we have $\s(x)=e'{}\i(X)\in X=e(x)$.
Thus $\s(x)\in e(x)$ for any $x\in V$. From (a) we have $x'\le x$ for any $x'\in e(x)$. Hence
$\s(x)\le x$ for any $x\in V$.
In a finite partially ordered set $Z$ any bijection $a:Z@>>>Z$ such that $a(z)\le z$ for all $z$ must be the
identity map. It follows that $\s=1$ so that $e=e'$. This proves (b).

\subhead 1.23\endsubhead
In 1.24-1.26 we describe the bijection in 1.17(c) assuming that $D$ is $2$, $4$ or $6$.
In each case we give a table in which there is one row for each $B\in\SS_D$; the row corresponding to $B$
is of the form $<B>:(\do)$ where $B$ is represented by the list of intervals of $B$ (we write an interval
such as $[4,6]$ as $456$) and $(\do)$ is a list of the vectors in $\la B\ra$ (we write
$1235$ instead of $e_1+e_2+e_3+e_5$, etc). In each list $(\do)$ we single out the vector $\e(B)$
in 1.17(b) by putting it in a box. Any non-boxed entry in $(\do)$ appears as a boxed entry in some previous row.
These tables extend the tables in \cite{\NEW, 1.17}.

\mpb

\subhead 1.24. The table for $D=2$\endsubhead

$\emp:(\bx{0})$

$<1>:(0,\bx{1})$

$<2>:(0,\bx{2})$.

$<12>:(0,\bx{12})$.

\subhead 1.25. The table for $D=4$\endsubhead

$\emp:(\bx{0})$

$<1>:(0,\bx{1})$

$<2>:(0,\bx{2})$

$<3>:(0,\bx{3})$

$<4>:(0,\bx{4})$

$<1,3>:(0,1,3,\bx{13})$

$<1,4>:(0,1,4,\bx{14})$

$<2,4>:(0,2,4,\bx{24})$ 

$<2,123>:(0,2,13,\bx{123})$

$<3,234>:(0,3,24,\bx{234})$

$<1234>:(0,\bx{1234})$

$<3,1234>:(0,3,1234,\bx{124})$

$<2,1234>:(0,2,1234,\bx{134})$

$<4,12>:(0,4,124,\bx{12})$

$<1,34>:(0,1,134,\bx{34})$

$<1234,23>:(0,1234,14,\bx{23})$.

\subhead 1.26. The table for $D=6$\endsubhead

$\emp:(\bx{0})$

$<1>:(0,\bx{1})$

$<2>:(0,\bx{2})$

$<3>:(0,\bx{3})$

$<4>:(0,\bx{4})$

$<5>:(0,\bx{5})$

$<6>:(0,\bx{6})$

$<1,4>:(0,1,4,\bx{14})$

$<1,6>:(0,1,6,\bx{16})$

$<2,4>:(0,2,4,\bx{24})$

$<2,5>:(0,2,5,\bx{25})$

$<2,6>:(0,2,6,\bx{26})$

$<3,6>:(0,3,6,\bx{36})$

$<4,6>:(0,4,6,\bx{46})$

$<1,3>:(0,1,3,\bx{13})$

$<1,5>:(0,1,5,\bx{15})$

$<3,5>:(0,3,5,\bx{35})$

$<2,123>:(0,2,13,\bx{123})$

$<3,234>:(0,3,24,\bx{234})$

$<4,345>:(0,4,35,\bx{345})$

$<5,456>:(0,5,46,\bx{456})$

$<1,3,5>:(0,1,3,5,13,15,35,\bx{135})$

$<1,3,6>:(0,1,3,6,13,16,36,\bx{136})$

$<1,4,345>:(0,1,4,345,14,35,135,\bx{1345})$

$<1,4,6>:(0,1,4,6,14,16,46,\bx{146})$

$<2,4,6>:(0,2,4,6,24,26,46,\bx{246})$

$<1,5,456>:(0,1,5,456,15,46,146,\bx{1456})$

$<2,5,456>:(0,2,5,456,25,46,246,\bx{2456})$

$<2,5,123>:(0,2,5,123,25,13,135,\bx{1235})$

$<2,6,123>:(0,2,6,123,26,13,136,\bx{1236})$

$<2,4,12345>:(0,2,4,24,1345,1235,135,\bx{12345})$    

$<3,234,12345>: (0,3,234,12345,24,15,135,\bx{1245})$

$<3,6,234>:(0,3,6,234,24,36,246,\bx{2346})$

$<3,5,23456>:(0,3,5,2456,35,2346,246,\bx{23456})$

$<4,345,23456>:(0,4,345,23456,35,26,246,\bx{2356})$.

$<123456>:(0,\bx{123456})$

$<5,123456>:(0,5,123456,\bx{12346})$

$<4,123456>:(0,4,123456,\bx{12356})$

$<3,123456>:(0,3,123456,\bx{12456})$

$<2,123456>:(0,2,123456,\bx{13456})$

$<6,1234>:(0,6,12346,\bx{1234})$

$<1,3456>:(0,1,13456,\bx{3456})$

$<2,5,123456>:(0,2,5,25,123456,13456,12346,\bx{1346})$

$<3,5,123456>:(0,3,5,35,123456,12456,12346,\bx{1246})$

$<2,4,123456>:(0,2,4,24,123456,13456,12356,\bx{1356})$

$<3,6,1234>:(0,3,6,36,1234,12346,1246,\bx{124})$

$<1,4,3456>:(0,1,4,14,3456,13456,1356,\bx{356})$

$<2,6,1234>:(0,2,6,26,1234,12346,1346,\bx{134})$

$<1,5,3456>:(0,1,5,14,3456,13456,1346,\bx{346})$

$<3,234,123456>:(0,3,234,24,123456,12456,1356,\bx{156})$

$<4,345,123456>:(0,4,345,35,123456,12356,1246,\bx{126})$

$<4,6,12>:(0,4,6,46,124,126,1246,\bx{12})$

$<1,3,56>:(0,1,3,13,156,356,1356,\bx{56})$

$<1,6,34>:(0,1,6,16,134,346,1346,\bx{34})$

$<5,12,456>:(0,5,12,456,12456,46,1246,\bx{125})$

$<2,56,123>:(0,2,56,123,12356,13,1356,\bx{256})$

$<123456,2345>:(0,123456,16,\bx{2345})$

$<123456,3,2345>:(0,123456,3,2345,12456,16,136,\bx{245})$

$<123456,4,2345>:(0,123456,4,2345,12356,16,146,\bx{235})$

$<123456,2,45>:(0,2,123456,13456,1236,245,136,\bx{45})$

$<123456,5,23>:(0,5,123456,12346,1456,235,146,\bx{23})$

$<3456,1,45>:(0,1,45,3456,13456,36,136,\bx{145})$

$<1234,6,23>:(0,6,23,1234,12346,14,146,\bx{236})$

$<123456,2345,34>:(0,123456,2345,34,16,25,1346,\bx{1256})$.

\subhead 1.27\endsubhead
For $m\in\NN$ such that $m\le D/2$ let $\SS^m_D=\{B\in\SS_D;|B^0|=m\}$. One can show:

(a) {\it $|\SS^m_D|=\binom{D+1}{(D/2)-m}$.}
\nl
Indeed $\SS^m_D$ can be identified with a fibre of the map $\ti u:V@>>>\NN$ in 1.11 and that fibre is in
bijection with a set of symbols with fixed defect as in \cite{\SYMP}. These symbols can be counted and we find
(a).

If $B\in\SS_D$, then $B^1\in\SS^0_D$. 
This is seen by induction on $D$. Alternatively, $B^1$ satisfies $(P_0),(P_1),(P_2)$ 
hence is in $\SS_D$, by 1.3(c).
Thus $B\m B^1$ is a well defined (surjective) map $\SS_D@>>>\SS^0_D$. One can show:

(b) {\it This map induces a bijection $\{B\in\SS_D;|B|=D/2\}@>\si>>\SS^0_D$.}

\subhead 1.28\endsubhead
We now assume that $G$ in 0.1 is of type $B_n$ or $C_n$, $n\ge2$, or $D_n$, $n\ge4$.
We define the set $\ti\BB_c$ in 0.1.
If $|c|=1$, $\ti\BB_c$ consists of $(1,1)$.  
Assume now that $|c|\ge2$. We associate to $c$ a number $D\in2\NN$, and an $\FF_2$-vector space $V$ 
with basis $\{e_i;i\in[1,D]\}$ as in 1.9 so that 
$\Irr_c$ is identified with $M(\cg_c)=V$ as in \cite{\ORA}.
Then $\CC[M(\cg_c)]$ becomes the vector space of functions $V@>>>\CC$. The elements of $\ti\BB_c$ are the
characteristic functions of the subsets $\la B\ra$ of $V$ for various $B\in\SS_D$. This has the properties 
(I)-(V) in 0.1. (The bipositivity property (I)
in 0.1 follows from the fact that $\la B\ra$ is an isotropic subspace of $V$ for any
$B\in\SS_D$.)

\head 2. The case where $D$ is odd\endhead
\subhead 2.1\endsubhead
In this section we will sketch without proof a variant of the definitions and results in \S1 in which 
$D\in\NN$ is taken to be odd. 
 
We say that $B\in R_D$ is primitive if either $B=\emp$ or $B$ is of the form

(a) $B=\{[1,D-1],[2,D-2],[k,D-k]\}$ for some odd $k\in\NN$ such that  $k\le(D-1)/2$.
\nl
We define a subset $\SS_D$ of $R_D$ by induction on $D$ as follows.
If $D=1$, $\SS_D$ consists of a single element namely $\emp\in R_D$.
If $D\ge3$  we say that $B\in R_D$ is in $\SS_D$ if either $B$ is primitive, or

(b) $|B^0|\ne0$ and there exists $i\in[1,D]$ and $B'\in\SS_{D-2}$ such that $B=t_i(B')$, or

(c) $|B^0|=0$ and there exists $i\in[1,D-1]$ and $B'\in\SS_{D-2}$ such that $B=t_i(B')$.

Here $t_i$ is as in 1.1.

\subhead 2.2\endsubhead
We shall use the notation of 1.9 (with $D$ odd). 
Let $\un V=V/\FF_2\z$ where $\z=e_1+e_3+e_5+\do+e_D$. Now $(,):V\T V@>>>\FF_2$ induces a nondegenerate
symplectic form $\un V\T\un V@>>>\FF_2$. Let $\p:V@>>>\un V$ be the obvious map. Now 
$\un V$ with its basis $\{\p(e_i);i\in[0,D-1]\}$ is like $V$ in 1.9 (of even dimension). Hence
$\pmb\cf(\un V)$ is defined and we have canonical bijections $\un\a:\SS_{D-1}@>\si>>\pmb\cf(\un V)$
(as in 1.17(c)) and $\un e:\un V@>\si>>\pmb\cf(\un V)$ (as in 1.22(b)).

For $B\in\SS_D$ let $\la B\ra$ be the subspace of $V$
generated by $\{e_I;I\in B\}$; this is in fact a basis of $\la B\ra$ and
$\{\p(e_I);I\in B\}$ is a basis of $\p(\la B\ra)$.
Let $\pmb\cf(V)$ be the set of (isotropic) subspaces of $\un V$ of the form $\p(\la B\ra)$ for some $B\in\SS_D$.
Now $\pmb\cf(V)$ does not in general coincide with $\pmb\cf(\un V)$.

On can show that the map $\a:\SS_D@>>>\pmb\cf(V)$, $B\m\p(\la B\ra)$ is a bijection and that
there is a unique bijection $e:\un V@>\si>>\pmb\cf(V)$ such that for any $x\in\un V$ we have $x\in e(x)$.
Consider the matrix indexed by $\un V\T\un V$ whose entry at $(x,x')\in\un V\T\un V$ is $1$ if $x'\in e(x)$
and is $0$ if $x'\n e(x)$. One can show that this matrix is upper triangular with $1$ on diagonal for a
suitable partial order on $\un V$.

\subhead 2.3\endsubhead
For $m\in\NN$ such that $m\le(D-1)/2$ let $\SS^m_D=\{B\in\SS_D;|B^0|=m\}$. 
For $m>0$, even, we have $\SS^m_D=\emp$.
One can show that the bijection $\a\i e\un e\i\un\a:\SS_{D-1}@>\si>>\SS_D$ (see 2.2)
restricts to the identity map $\SS_D^0@>>>\SS_{D-1}^0$ and to a bijection
$\SS^m_{D-1}\cup\SS^{m+1}_{D-1}@>\si>>\SS^m_D$ for $m$ odd. 

\subhead 2.4\endsubhead
In 2.5-2.7 we describe the bijection $\SS_D@>\si>>\pmb\cf(V)$, $B\m\p(\la B\ra)$ in 2.2 assuming 
that $D$ is $3,5$ or 
$7$. In each case we give a table in which there is one row for each $B\in\SS_D$; the row corresponding to $B$
is of the form $<B>:(\do)$ where $B$ is represented by the list of intervals of $B$. We use conventions
similar to those in 1.23, except that now $(\do)$ is a list of vectors in $\un V$
(we write $1235$ instead of $\p(e_1)+\p(e_2)+\p(e_3)+\p(e_5)$, etc). In each list $(\do)$ we single out 
(by putting it in a box) the vector $x\in\un V$ such that $e(x)=\p(\la B\ra)$ with $e$ as in 2.2.
Any non-boxed entry in $(\do)$ appears as a boxed entry in some previous row.

\subhead 2.5. The table for $D=3$\endsubhead

$\emp:(\bx{0})$

$<1>:(0,\bx{1})$

$<2>:(0,\bx{2})$.

$<12>:(0,\bx{12})$.

\subhead 2.6. The table for $D=5$\endsubhead

$\emp:(\bx{0})$

$<1>:(0,\bx{1})$

$<2>:(0,\bx{2})$

$<3>:(0,\bx{3})$

$<4>:(0,\bx{4})$

$<1,3>:(0,1,3,\bx{13})$

$<1,4>:(0,1,4,\bx{14})$

$<2,4>:(0,2,4,\bx{24})$ 

$<2,123>:(0,2,13,\bx{123})$

$<3,234>:(0,3,24,\bx{234})$

$<1234>:(0,\bx{1234})$

$<3,1234>:(0,3,1234,\bx{124})$

$<2,1234>:(0,2,1234,\bx{134})$

$<4,12>:(0,4,124,\bx{12})$

$<1,34>:(0,1,134,\bx{34})$

$<5,12>:(0,12,13,\bx{23})$.

\subhead 2.7. The table for $D=7$\endsubhead

$\emp:(\bx{0})$

$<1>:(0,\bx{1})$

$<2>:(0,\bx{2})$

$<3>:(0,\bx{3})$

$<4>:(0,\bx{4})$

$<5>:(0,\bx{5})$

$<6>:(0,\bx{6})$

$<1,4>:(0,1,4,\bx{14})$

$<1,6>:(0,1,6,\bx{16})$

$<2,4>:(0,2,4,\bx{24})$

$<2,5>:(0,2,5,\bx{25})$

$<2,6>:(0,2,6,\bx{26})$

$<3,6>:(0,3,6,\bx{36})$

$<4,6>:(0,4,6,\bx{46})$

$<1,3>:(0,1,3,\bx{13})$

$<1,5>:(0,1,5,\bx{15})$

$<3,5>:(0,3,5,\bx{35})$

$<2,123>:(0,2,13,\bx{123})$

$<3,234>:(0,3,24,\bx{234})$

$<4,345>:(0,4,35,\bx{345})$

$<5,456>:(0,5,46,\bx{456})$

$<1,3,5>:(0,1,3,5,13,15,35,\bx{135})$

$<1,3,6>:(0,1,3,6,13,16,36,\bx{136})$

$<1,4,345>:(0,1,4,345,14,35,135,\bx{1345})$

$<1,4,6>:(0,1,4,6,14,16,46,\bx{146})$

$<2,4,6>:(0,2,4,6,24,26,46,\bx{246})$

$<1,5,456>:(0,1,5,456,15,46,146,\bx{1456})$

$<2,5,456>:(0,2,5,456,25,46,246,\bx{2456})$

$<2,5,123>:(0,2,5,123,25,13,135,\bx{1235})$

$<2,6,123>:(0,2,6,123,26,13,136,\bx{1236})$

$<2,4,12345>:(0,2,4,24,1345,1235,135,\bx{12345})$    

$<3,234,12345>: (0,3,234,12345,24,15,135,\bx{1245})$

$<3,6,234>:(0,3,6,234,24,36,246,\bx{2346})$

$<3,5,23456>:(0,3,5,2456,35,2346,246,\bx{23456})$

$<4,345,23456>:(0,4,345,23456,35,26,246,\bx{2356})$

$<123456>:(0,\bx{123456})$

$<5,123456>:(0,5,123456,\bx{12346})$

$<4,123456>:(0,4,123456,\bx{12356})$

$<3,123456>:(0,3,123456,\bx{12456})$

$<2,123456>:(0,2,123456,\bx{13456})$

$<6,1234>:(0,6,12346,\bx{1234})$

$<1,3456>:(0,1,13456,\bx{3456})$

$<2,5,123456>:(0,2,5,25,123456,13456,12346,\bx{1346})$  

$<3,5,123456>:(0,3,5,35,123456,12456,12346,\bx{1246})$  

$<2,4,123456>:(0,2,4,24,123456,13456,12356,\bx{1356})$  

$<3,6,1234>:(0,3,6,36,1234,12346,1246,\bx{124})$    

$<1,4,3456>:(0,1,4,14,3456,13456,1356,\bx{356})$    

$<2,6,1234>:(0,2,6,26,1234,12346,1346,\bx{134})$   

$<1,5,3456>:(0,1,5,14,3456,13456,1346,\bx{346})$   

$<3,234,123456>:(0,3,234,24,123456,12456,1356,\bx{156})$ 

$<4,345,123456>:(0,4,345,35,123456,12356,1246,\bx{126})$ 

$<4,6,12>:(0,4,6,46,124,126,1246,\bx{12})$     

$<1,3,56>:(0,1,3,13,156,356,1356,\bx{56})$    

$<1,6,34>:(0,1,6,16,134,346,1346,\bx{34})$    

$<5,12,456>:(0,5,12,456,12456,46,1246,\bx{125})$

$<2,56,123>:(0,2,56,123,12356,13,1356,\bx{256})$

 $<7,1234>:(0,135,1234,\bx{245})$ 

$<3,7,1234>:(0,3,15,1234,135,124,245,\bx{2345})$

$<4,7,12>:(0,4,135,12,124,1345,2345,\bx{235})$

 $<2,7,1234>:(0,2,135,1234,1235,134,245,\bx{45})$

 $<5,7,12>:(0,5,13,12,135,125,235,\bx{23})$

$<1,7,34>:(0,1,34,135,134,35,45,\bx{145})$

$<6,12,567>:(0,6,12,136,13,126,23,\bx{236})$

$<123456,2345,34>:(0,123456,2345,34,16,25,1346,\bx{1256})$.

\head 3. Exceptional groups\endhead
\subhead 3.1\endsubhead
Let $\G$ be a finite group. Let $x\in\G$ and let $\r$ be a not necessarily irreducible representation over
$\CC$ of the centralizer $Z_\G(x)$ of $x$ in $\G$. We define $(x,\r)\in M(\G)$ to be 
$\sum_\s(\s:\r)(x,\s)$ where $\s$ runs over the irreducible representations of $Z_\G(x)$
up to isomorphism and $:$ denotes multiplicity. 
Let $H$ be a subgroup of $\G$. Following \cite{\ORA, p.312} we define a linear map 
$i_{H,\G}:\CC[M(H)]@>>>\CC[M(\G)]$ by 
$$(x,\s)\m(x,\Ind_{Z_H(x)}^{Z_\G(x)}(\s)).\tag a$$
As stated in {\it loc.cit.} we have  

(b) {\it $i_{H,\G}(A_H(f))=A_\G(i_{H,\G}(f))$ for any $f\in\CC[M(H)]$.}
\nl
If $f\in \CC[M(H)]$ is $\ge0$ then clearly $i_{H,\G}(f)$ is $\ge0$. Using this and (b) we see that

(c) {\it If $f\in \CC[M(H)]$ is bipositive then $i_{H,\G}(f)\in\CC[M(\G)]$ is bipositive.}
\nl
Assume now that $H$ is a normal subgroup of $\G$ and let $\p:\G@>>>\G/H$ be the canonical map.
Following {\it loc.cit.} we define a linear map $\p_{H,\G}:\CC[M(\G/H)]@>>>\CC[M(\G)]$ by 
$$(x,\s)\m\sum_{y\in\p\i(x)}\sum_{\t\in\Irr(Z_\G(y))}|Z_\G(y)||Z_{\G/H}(x)|\i|\G|\i|\G/H|(\t:\s)(y,\t)\tag d$$    
where $\t$ runs over the irreducible representations of $Z_\G(y)$ up to isomorphism and $\t:\s$ denotes the 
multiplicity of $\t$ in $\s$ viewed as a representation of $Z_\G(y)$ via the obvious homomorphism 
$Z_\G(y)@>>>Z_{\G/H}(x)$. 
As stated in {\it loc.cit.} we have  

(e) {\it $\p_{H,\G}(A_{\G/H}(f))=A_\G(\p_{H,\G}(f))$ for any $f\in\CC[M(\G/H)]$.}
\nl
If $f\in \CC[M(\G/H)]$ is $\ge0$ then clearly $\p_{H,\G}(f)$ is $\ge0$. Using this and (e) we see that

(f) {\it If $f\in \CC[M(\G/H)]$ is bipositive then $\p_{H,\G}(f)\in\CC[M(\G)]$ is bipositive.}

Now let $H\sub H'$ be two subgroups of $\G$ such that $H$ is normal in $H'$.
We define a linear map $\ss_{H,H'}:\CC[M(H'/H)]@>>>\CC[M(\G)]$ by $f\m i_{H',\G}(\p_{H,H'}(f))$.
From (c),(f) we deduce:

(g) {\it If $f\in\CC[M(H'/H)]$ is bipositive then $\ss_{H,H'}(f)\in\CC[M(\G)]$ is bipositive.}
\nl
Note that $\ss_{H,H'}(1,1)$ is the same as $S_{H,H'}$ defined in \cite{\NEW}; in this special
case (g) can be also deduced from \cite{\NEW, 0.7}.

\subhead 3.2\endsubhead
For $N\ge1$ let $S_N$ be the group of all permutations of $[1,N]$.
We shall use the notation of \cite{\ORA, 4.3} for the elements of $M(S_N)$ with $N=2,3,4$ or $5$
(but we replace $\bar Q_l$ by $\CC$).
We now give some examples of bipositive elements. Note that $(1,1)\in M(\G)$ is bipositive for any 
finite group $\G$. Indeed, we have
$$A_\G(1,1)=\sum_{(x,\s)\in M(\G)}\dim\s|Z_\G(z)|\i(x,\s).$$
Let

$\L_{-1}=(g_2,\e)+(1,1)\in M(S_2)$

$\L'_{\th^j}=(g_3,\th^j)+(g_2,1)+(1,1)\in M(S_3)$, ($j=1,2$)

$\L_{\th^j}=(g_3,\th^j)+(g_2,\e)+(1,1)\in M(S_3)$, ($j=1,2$)

$\L_{i^k}=(g_4,i^k)+(g_4,-1)+(g_3,1)+(1,\l^2)+(1,1)\in M(S_4)$, ($k=1,-1$).

$\L_{\z^j}=(g_5,\z^j)+(1,\l^4)+2(1,\l^2)+(1,\nu)+(1,\nu')+(1,1)\in M(S_5)$,  ($j=1,2,3,4$).

$$\align&\L'_{\z^l,\z^{2l}}=(g_5,\z^l)+(g_5,\z^{2l})
+(g'_2,1)+(g'_2,\e')\\&+(g'_2,\e'')+(g'_2,\e)+(1,\l^2)+(1,\nu)+(1,1)\in M(S_5),\qua l=1,2,3,4.
\endalign$$
Here $\th=\exp(2\pi i/3),\z=\exp(2\pi i/5)$.

One can verify by computation that each of the elements above 
(except for $\L_{\z^j}$) is fixed by the non-abelian Fourier
transform hence is bipositive. In 3.3 we will show that $\L_{\z^j}$ is also bipositive.
We say that 

$(1,1)$ is  the primitive element of $M(S_1)$;

$\L_{-1},(1,1)$ are the primitive elements of $M(S_2)$;

$\L'_{\th^j},(1,1)$ are the primitive elements of $M(S_3)$ (when $G$ is not simply laced);

$\L_{\th^j},(1,1)$ are the primitive elements of $M(S_3)$ (when $G$ is simply laced);

$\L_{i^k},(1,1)$ are the primitive elements of $M(S_4)$;

$\L_\z,\L'_{\z,\z^2},\L'_{\z^2,\z^4},\L'_{\z^3,\z},(1,1)$ are the primitive elements of $M(S_5)$.

It follows that the following elements are bipositive.

$\L_{-1,-1}=\L_{-1}\bxt\L_{-1}\in M(S_2)\ot M(S_2)=M(S_2\T S_2)$

$\L_{-1,1}=\L_{-1}\bxt (1,1)\in M(S_2)\ot M(S_2)=M(S_2\T S_2)$

$\L_{1,-1}=(1,1)\bxt \L_{-1}\in M(S_3)\ot M(S_2)=M(S_3\T S_2)$

$\L_{\th^j,-1}=\L_{\th^j}\bxt\L_{-1}\in M(S_3)\ot M(S_2)=M(S_3\T S_2)$, ($j=1,2$)

$\L_{\th^j,1}=\L_{\th^j}\bxt (1,1)\in M(S_3)\ot M(S_2)=M(S_3\T S_2)$, ($j=1,2$)

Note that both $\L_{-1,-1},\L_{\th^j,-1}$ are fixed by the non-abelian Fourier transform. 
We say that 

$\L_{-1,-1},\L_{-1,1},(1,1)$ are the primitive elements of $M(S_2\T S_2)$;

$\L_{\th^j,-1},\L_{\th^j,1},\L_{1,-1},(1,1)$ are the primitive elements of $M(S_3\T S_2)$.

\subhead 3.3\endsubhead
Let $H$ be a dihedral group of order $10$. We denote by $g_5$ an element of order $5$ of $H$ and by
$g_2$ an element of order $2$ such that $g_2g_5g_2\i=g_5\i$. Now $H$ has four conjugacy classes; they have
 representatives $1,g_2,g_5,g_5^2$ with centralizers of order $10,2,5,5$.
The irreducible representations of $H$ are $1,r,r',\e$ where $r,r'$ are $2$-dimensional and $\e$ is the sign.
We can assume that $\tr(g_5,r)=\tr(g_5,r')=\z+\z\i$, $\tr(g_5^2,r)=\tr(g_5^2,r')=\z^2+\z^{-2}$,
$\tr(g_2,r)=\tr(g_2,r')=0$, $\tr(g_5,\e)=\tr(g_5^2,\e)=1$, $\tr(g_2,\e)=-1$. The elements of $M(H)$ are
$(1,1),(1,r),(1,r'),(1,\e),(g_2,1),(g_2,\e)$, $(g_5^k,\z^l)$ with $k=1,2$, $l=0,1,\do,4$.
Here $\z^l$ is the character of the cyclic group generated by $g_5$ which takes the value $\z^l$ at $g_5$.
For $C\in\ZZ$ we set $[C]=\z^C+\z^{-C}$. Note that $[C]$ depends only on the residue class of $C$ modulo $5$. 
We write $A$ instead of $A_H$. We have
$$\align&A(1,1)=(1/10)(1,1)+(1/5)(1,r)+(1/5)(1,r')+(1/10)(1,\e)+(1/2)(g_2,1)\\&+(1/2)(g_2,\e)+
\sum_{k'\in\{1,2\},l'\in\{0,4\}}(1/5)(g^{k'},\z^{l'})\endalign$$
$$\align&A(1,\e)=(1/10)(1,1)+(1/5)(1,r)+(1/5)(1,r')+(1/10)(1,\e)-(1/2)(g_2,1)\\&-(1/2)(g_2,\e)+
\sum_{k'\in\{1,2\},l'\in\{0,4\}}(1/5)(g^{k'},\z^{l'})\endalign$$
$$A(g_2,1)=(1/2)(1,1)-(1/2)(1,\e)+(1/2)(g_2,1)-(1/2)(g_2,\e)$$
$$\align&A(g_5^k,\z^l)=(1/5)(1,1)+(1/5)[k](1,r)+(1/5)[2k](1,r')+(1/5)(1,\e)\\&+\sum_{k'\in\{1,2\},l'\in\{0,4\}}
(1/5)[kl'-k'l](g^{k'},\z^{l'}).\endalign$$
Assume that $k=1$ and $l\in[1,4]$. Using $[1]+[2]=-1$, $[2]+[4]=-1$, we have
$$\align&A(g_5,\z^l)+A(g_5^2,\z^{2l})
=(2/5)(1,1)-(1/5)(1,r)-(1/5)(1,r')+(2/5)(1,\e)\\&+\sum_{k'\in\{1,2\},l'\in\{0,4\}}(1/5)([l'-k'l]+[2l'-2k'l])
(g^{k'},\z^{l'}).\endalign$$
Let $N_1=l'-k'l$, $N_2=2N_1$. If $N_1=0\mod5$ we have $[N_1]+[N_2]=[0]+[0]=4$.
Assume now that $N_1\ne0\mod5$. If $N_1+N_2=0\mod5$ then $3N_1=0\mod5$ so that $N_1=0\mod5$ contradicting our 
assumption. Thus $N_1,N_2$ are $\ne0$ in $\ZZ/5$ and their sum is $\ne0$ in $\ZZ/5$. This implies that
$[N_1]+[N_2]=[1]+[2]=-1$. We see that
$$\align&A(g_5,\z^l)+A(g_5^2,\z^{2l})
=(2/5)(1,1)-(1/5)(1,r)-(1/5)(1,r')+(2/5)(1,\e)\\&+
(4/5)(g_5,\z^l)+(4/5)(g_5^2,\z^{2l})+
+\sum_{k'\in\{1,2\},l'\in\{0,4\};l'-k'l\ne0\mod5}(-1/5)(g^{k'},\z^{l'}).\endalign$$
Hence
$$\align&A(g_5,\z^l)+A(g_5^2,\z^{2l})+A(g_2,1)+A(1,1)\\&
=(2/5)(1,1)-(1/5)(1,r)-(1/5)(1,r')+(2/5)(1,\e)\\&+
(4/5)(g_5,\z^l)+(4/5)(g_5^2,\z^{2l})
+\sum_{k'\in\{1,2\},l'\in\{0,4\};l'-k'l\ne0\mod5}(-1/5)(g^{k'},\z^{l'})\\&
+(1/2)(1,1)-(1/2)(1,\e)+(1/2)(g_2,1)-(1/2)(g_2,\e)\\&
+(1/10)(1,1)+(1/5)(1,r)+(1/5)(1,r')+(1/10)(1,\e)+(1/2)(g_2,1)+(1/2)(g_2,\e)\\&+
\sum_{k'\in\{1,2\},l'\in\{0,4\}}(1/5)(g^{k'},\z^{l'})=
(g_5,\z^l)+(g_5^2,\z^{2l})+(g_2,1)+(1,1)\endalign$$
that is
$$(g_5,\z^l)+(g_5^2,\z^{2l})+(g_2,1)+(1,1)\text{ is fixed by }A.$$

Next we show that the coefficient of any basis element $(x,\s)$ in
$$\align&A(g_5^k,\z^l)+A(1,\e)+A(1,1)=
(1/5)(1,1)+(1/5)[k](1,r)+(1/5)[2k](1,r')\\&+(1/5)(1,\e)+\sum_{k'\in\{1,2\},l'\in\{0,4\}}
(1/5)[kl'-k'l](g^{k'},\z^{l'})\\&
(1/10)(1,1)+(1/5)(1,r)+(1/5)(1,r')+(1/10)(1,\e)+(1/2)(g_2,1)+(1/2)(g_2,\e)\\&+
\sum_{k'\in\{1,2\},l'\in\{0,4\}}(1/5)(g^{k'},\z^{l'})+(1/10)(1,1)+(1/5)(1,r)+(1/5)(1,r')
\\&+(1/10)(1,\e)-(1/2)(g_2,1)-(1/2)(g_2,\e)+
\sum_{k'\in\{1,2\},l'\in\{0,4\}}(1/5)(g^{k'},\z^{l'})\endalign$$
is $\ge0$.
It is enough to show that if $k'\in\{1,2\},l'\in\{0,4\}$ then
$[kl'-k'l]+2\ge0$ and that $[k]+2\ge0$, $[2k]+2\ge0$. More generally, for any $C\in\ZZ$ we have $[C]+2\ge0$.

\mpb

We can regard $H$ as a subgroup of $S_5$ so that $g_5\in H$ becomes a $5$-cycle $g_5\in S_5$.
Then $\ss_{1,H}:M(H)@>>>M(S_5)$ is defined and for $l\in[1,4]$ we have
$$\ss_{1,H}((g_5,\z^l)+(g_5^2,\z^{2l})+(g_2,1)+(1,1))=\L'_{\z^l,\z^{2l}}\in M(S_5),\tag a$$
$$\ss_{1,H}((g_5,\z^l)+(1,\e)+(1,1))=\L_{\z^l}\in M(S_5).\tag b$$
It follows that the elements (a),(b) are bipositive. (The element (a) is fixed by $A_{S_5}$.)

\subhead 3.4\endsubhead
In the remainder of this section we assume that $G$ in 0.1 is of exceptional type. We are in one of the 
following cases:

(i) $|c|=1$, $\cg_c=S_1$.

(ii) $|c|=2$ (with $W$ of type $E_7$ or $E_8$), $\cg_c=S_2$.

(iii) $|c|=3$, $\cg_c=S_2$.

(iv) $|c|=4$ (with $W$ of type $G_2$), $\cg_c=S_3$.

(v) $|c|=5$ (with $W$ of type $E_6.E_7$ or $E_8$), $\cg_c=S_3$.

(vi) $|c|=11$ (with $W$ of type $F_4$), $\cg_c=S_4$.

(vii) $|c|=17$ (with $W$ of type $E_8$), $\cg_c=S_5$.

\subhead 3.5\endsubhead
In the case 3.4(i) we define $\ti\BB_c$ as the set consisting of $(1,1)\in M(S_1)$. 

In the cases 3.4(ii),3.4(iii) we define $\ti\BB_c$ as the subset of $\CC[M(S_2)]$ consisting of 

$\widehat{(1,1)}=\ss_{1,S_2}(1,1)=(1,1)$,        

$\widehat{(g_2,1)}=\ss_{S_2,S_2}(1,1)=(g_2,1)+(1,1)$,

$\widehat{(1,\e)}=\ss_{1,1}(1,1)=(1,\e)+(1,1)$,

$\widehat{(g_2,\e)}=\ss_{1,S_2}(\L_{-1})=\L_{-1}=(g_2,\e)+(1,1)$.

\subhead 3.6\endsubhead
In cases 3.4(iv),(v) we define $\ti\BB_c$ as the subset of $\CC[M(S_3)]$  consisting of 

$\widehat{(1,1)}=\ss_{1,S_3}(1,1)=(1,1)$,

$\widehat{(1,r)}=\ss_{1,H_{21}}(1,1)=(1,r)+(1,1)$,

$\widehat{(g_2,1)}=\ss_{H_{21},H_{21}}(1,1)=(g_2,1)+(1,r)+(1,1)$,

$\widehat{(g_3,1)}=\ss_{S_3,S_3}(1,1)=(g_3,1)+(g_2,1)+(1,1)$,

$\widehat{(1,\e)}=\ss_{1,1}(1,1)=(1,\e)+2(1,r)+(1,1)$,

$\widehat{(g_2,\e)}=\ss_{1,H_{21}}\L_{-1}=(g_2,\e)+(1,r)+(1,1)$,

and of

$\widehat{(g_3,\th^j)}=\ss_{1,S_3}\L'_{\th^j}=(g_3,\th^j)+(g_2,1)+(1,1)$,  $(j=1,2)$ (in case 3.4(iv))

$\widehat{(g_3,\th^j)}=\ss_{1,S_3}\L_{\th^j}=(g_3,\th^j)+(g_2,\e)+(1,1)$,  $(j=1,2)$ (in case 3.4(v)).

Here the index $H,H'$ in $\ss_{H,H'}$ is a pair of subgroups of $S_3$ as in \cite{\NEW,3.10}.

\subhead 3.7\endsubhead
In the case 3.4(vi) we define $\ti\BB_c$ as the subset of $\CC[M(S_4)]$  consisting of 

$\widehat{(1,1)}=\ss_{1,S_4}(1,1)$

$\widehat{(1,\l^1)}=\ss_{1,H_{31}}(1,1)$

$\widehat{(1,\s)}=\ss_{1,H_{22}}(1,1)$

$\widehat{(1,\l^2)}=\ss_{1,H_{211}}(1,1)$

$\widehat{(g_2,1)}=\ss_{\ti H_{211},H_{22}}(1,1)$

$\widehat{(g'_2,1)}=\ss_{\ti H_{22},\ti H}(1,1)$

$\widehat{(g_2,\e'')}=\ss_{H_{211},H_{221}}(1,1)$

$\widehat{(g_3,1)}=\ss_{H_{31},H_{31}}(1,1)$

$\widehat{(g_4,1)}=\ss_{S_4,S_4}(1,1)$

$\widehat{(g'_2,\e'')}=\ss_{H_{22},H_{22}}(1,1)$

$\widehat{(g'_2,\e')}=\ss_{\ti H,\ti H}(1,1)$   

$\widehat{(g_2,\e')}=\ss_{1,H_{22}}\L_{-1,1}=(g_2,\e')+(1,\s)+(1,\l^1)+(1,1)$

$\widehat{(g'_2,r)}=\ss_{\ti H_{211},H_{22}}\L_{-1}=(g_2,\e')+(g'_2,r)+(g_2,1)+(1,\l^1)+(1,\s)+(1,1)$

 $\widehat{(g_4,-1)}=\ss_{H_{22},\ti H}\L_{-1}=(g_4,-1)+(g'_2,r)+(g'_2,1)+(g_2,1)+(1,\s)+(1,1)$                        
$\widehat{(1,\l^3)}=\ss_{1,1}(1,1)=(1,\l^3)+3(1,\l^2)+3(1,\l^1)+2(1,\s)+(1,1)$
                         
$\widehat{(g_2,\e)}=\ss_{1,H_{211}}\L_{-1}=(g_2,\e)+(g_2,\e')+2(1,\l^1)+(1,\l^2)+(1,\s)+(1,1)$
                        
$\widehat{(g'_2,\e)}=\ss_{1,H_{22}}\L_{-1,-1}=(g'_2,\e)+(g'_2,1)+(g_2,\e')+(g_2,\e'')+(1,\l^1)+(1,\s)+(1,1)$
                         
$\widehat{(g_3,\th^j)}=\ss_{1,H_{31}}\L'_{\th^j}=(g_3,\th^j)+(g_2,1)+(g_2,\e')+(1,\l^1)+(1,1)$ ($j=1,2$)
                          
$\widehat{(g_4,i^k)}=\ss_{1,S_4}\L_{i^k}=(g_4,i)+(g_4,-1)+(g_3,1)+(1,\l^2)+(1,1)$, ($k=1,-1$).
                          
Here the index $H,H'$ in $\ss_{H,H'}$ is a pair of subgroups of $S_4$ as in \cite{\NEW,3.10} except that
$\ss_{1,1}$ does not appear there. In each case $H/H'$ is a product of symmetric groups.

Consider the  matrix (from \cite{\NEW}):
$$\left(\matrix
1  &0    &0   &0    &0   &0      &0   &0    &0  &0  &0                   \\
  1   &1  &0   &0    &0    &0     &0   &0    &0  &0   &0                 \\
 1   &1    &1  &0    &0    &0     &0   &0    &0  &0   &0               \\
 1    & 2  & 1   &1  & 0    &0    & 0   &0   & 0  & 0   &0                    \\
 1    & 1  & 1   &0   &1     &0   &0   &0    &0  &0   &0             \\
  1   & 0   &1   &0   & 1   &1    &0   &0    &0  &0   &0               \\
 1    & 2   &1   &1    &1    &0   &1   &0    &0  &0   &0               \\
 1    & 1   &0   & 0    &1    &0   & 1   &1   &0  &0   &0                \\
 1    & 0   &0   & 0    &1    &1   & 0    &1   &1 &0     &0            \\
 1    & 1   &1   & 0    &2    &1    &0    &0   & 0  &1   &0              \\
 1    & 0   &1    &0    &1    &2    &0    &0   & 1  &0   &1           \\
\endmatrix\right)$$
with rows indexed from left to right and columns indexed from up to down
by the elements of $M_0(S_4)$  in the order
$$(1,1),(1,\l^1),(1,\s),(1,\l^2),(g_2,1),(g'_2,1),(g_2,\e''),(g_3,1),(g_4,1),(g'_2,\e''),(g'_2,\e').$$

For $(x,\s)\in M_0(S_4)$, the coefficient of $(x',\s')\in M_0(S_4)$ in $\widehat{(x,\s)}\in\CC[M(S_4)]$ is 
 the entry of the matrix above in the row $(x,\s)$ and column $(x',\s')$;
the coefficient of any $(x',\s')\in M(S_4)-M_0(S_4)$ is $0$.

\subhead 3.8\endsubhead
In the case 3.4(vii) we define $\ti\BB_c$ as the subset of $\CC[M(S_5)]$  consisting of 

$\widehat{(1,1)}=\ss_{1,S_5}(1,1)$   

$\widehat{(1,\l^1)}=\ss_{1,H_{41}}(1,1)$

$\widehat{(1,\nu)}=\ss_{1,H_{32}}(1,1)$

$\widehat{(1,\l^2)}=\ss_{1,H_{311}}(1,1)$

$\widehat{(1,\nu')}=\ss_{1,H_{221}}(1,1)$

$\widehat{(1,\l^3)}=\ss_{1,H_{2111}}(1,1)$

$\widehat{(g_2,1)}=\ss_{\ti H_{2111},H_{32}}(1,1)$

$\widehat{(g_2,r)}=\ss_{\ti H_{2111},H_{221}}(1,1)$

$\widehat{(g_3,1)}=\ss_{\ti H_{311},H_{32}}(1,1)$

$\widehat{(g'_2,1)}=\ss_{\ti H_{221},\ti H}(1,1)$

$\widehat{(g'_2,\e'')}=\ss_{H_{221},H_{221}}(1,1)$

$ \widehat{(g_6,1)}=\ss_{H_{32},H_{32}}(1,1)$

$\widehat{(g_2,\e)}=\ss_{H_{2111},H_{2111}}(1,1)$

$\widehat{(g_3,\e)}=\ss_{H_{311},H_{311}}(1,1)$

$ \widehat{(g_4,1)}=\ss_{H_{41},H_{41}}(1,1)$

$\widehat{(g_5,1)}=\ss_{S_5,S_5}(1,1)$

$\widehat{(g'_2,\e')}=\ss_{\ti H,\ti H}(1,1)$

 $\widehat{(g_2,-1)}=\ss_{1,H_{32}}\L_{1,-1}=(g_2,-1)+(1,\l^1)+(1,\nu)+(1,1)$

$\widehat{(g_2,-r)}=\ss_{1,H_{221}}\L_{-1,1}=(g_2,-r)+(g_2,-1)+(1,\l^2)+(1,\nu')+2(1,\nu)+2(1,\l^1)+(1,1)$

$\widehat{(g'_2,r)}=\ss_{\ti H_{2111},H_{221}}\L_{-1}=
(g'_2,r)+(g_2,-r)+(g_2,-1)+(g_2,1)+(g_2,r)+(1,\l^2)+(1,\nu')+2(1,\nu)+2(1,\l^1)+(1,1)$

$\widehat{(g_4,-1)}=\ss_{\ti H_{221},\ti H}\L_{-1}=
(g_4,-1)+(g'_2,r)+(g'_2,1)+(g_2,r)+(g_2,1)+(1,\l^1)+(1,\nu)+(1,\nu')+(1,1)$
                   
$$\align&\widehat{(g_6,-1)}=\ss_{\ti H_{311},H_{32}}\L_{-1}\\&=
(g_6,-1)+(g'_2,r)+(g_2,-1)+(g_3,1)+(g_2,1)+(g_2,r)+(1,\l^1)+(1,\nu)+(1,1)\endalign$$

$\widehat{(g_3,\th^j)}=\ss_{1,H_{32}}\L_{1,\th^j}
=(g_3,\th^j)+(g_2,r)+(g_2,\e)+(1,\l^1)+(1,\nu)+(1,1)$, ($j=1,2$)

$$\align&\widehat{(g_6,\th^j)}=\ss_{\ti H_{2111},H_{32}}\L_{\th^j}=
(g_6,\th^j)\\&+(g_3,\th)+(g'_2,r)+(g_2,r)+(g_2,\e)+(g_2,1)+(1,\l^1)+(1,\nu)+(1,1),  (j=1,2)\endalign$$ 

 $\widehat{(1,\l^4)}=\ss_{1,1}(1,1)=4(1,\l^1)+6(1,\l^2)+4(1,\l^3)+(1,\l^4)+5(1,\nu)+5(1,\nu')+(1,1)$

$\widehat{(g_2,-\e)}=\ss_{1,H_{2111}}\L_{-1}=(g_2,-\e)+2(g_2,-r)+(g_2,-1)+3(1,\l^1)+3(1,\l^2)+(1,\l^3)+3(1,\nu)+2(1,\nu')+(1,1)$

$\widehat{(g_3,\e\th^j)}=\ss_{1,H_{311}}\L_{\th^j}
=(g_3,\e\th^j)+(g_3,\th)+(g_2,1)+2(g_2,r)+(g_2,\e)+2(1,\l^1)+(1,\l^2)+(1,\nu)+(1,1)$,
 ($j=1,2$)     

$\widehat{(g'_2,\e)}=\ss_{1,H_{221}}\L_{-1,-1}=(g'_2,\e)+(g'_2,1)+2(g_2,-1)+2(g_2,-r)+(1,\l^2)+(1,\nu')+2(1,\nu)+2(1,\l^1)+(1,1)$

$\widehat{(g_6,-\th^j)}=\ss_{1,H_{32}}\L_{\th^j,-1}
=(g_6,-\th^j)+(g_3,\th)+(g'_2,r)+(g_2,1)+(g_2,r)+(g_2,-1)+(1,\l^1)+(1,\nu)+(1,1)$,
 ($j=1,2$)

$\widehat{(g_4,i^k)}=\ss_{1,H_{41}}\L_{i^k}=
(g_4,i^k)+(g_4,-1)+(g_3,1)+(g_3,\e)+(1,\l^2)+(1,\l^3)+(1,\l^1)+(1,\nu)+(1,1)$, ($k=1,-1$)

$\widehat{(g_5,\z)}=\L_\z$,

$\widehat{(g_5,\z^2)}=\L'_{\z,\z^2}$,

$\widehat{(g_5,\z^3)}=\L'_{\z^3,\z}$,

$\widehat{(g_5,\z^4)}=\L'_{\z^2,\z^4}$.
\nl
Here the index $H,H'$ in $\ss_{H,H'}$ is a pair of subgroups of $S_5$ as in \cite{\NEW,3.10} except that
$\ss_{1,1}$ does not appear there. In each case $H/H'$ is a product of symmetric groups.

Consider the matrix (from \cite{\NEW}):
$$\left(\matrix
1  &0   &0    &0   &0   &0   &0   &0      &0  &0    &0   &0   &0   &0    &0   &0                 &0        \\
 1    &1 &0    &0   &0   &0   &0  &0     &0  &0     &0  &0  &0   &0    &0   &0            &0                  \\
 1   &1    &1  &0   &0   &0   &0   &0      &0  &0    &0  &0  &0   &0   &0   &0              &0                \\
 1    &2   &1   &1  &0   &0   &0   &0      &0  &0     &0 &0   &0   &0   &0   &0              &0               \\
 1    &2   &2    &1  &1  &0   &0   &0      &0  &0     &0 &0   &0   &0    &0   &0     &0               \\
 1    &3   &3    &3   &2  &1  &0   &0      &0  &0     &0  &0  &0   &0    &0   &0     &0               \\
 1    &1   &1    &0   &0   &0  &1  &0      &0  &0     &0  &0   &0   &0   &0   &0       &0            \\
 1    &2   &2    &1   &1   &0   &1  &1     &0  &0     &0  &0   &0   &0   &0   &0              &0          \\
 1    &1   &1    &0   &0   &0   &1   &1     &1&0      &0  &0   &0   &0    &0   &0               &0        \\
 1    &1   &1    &0   &1   &0   &1   &1      &0  &1   &0  &0   &0   &0    &0   &0              &0            \\
 1    &2   &2    &1   &1   &0   &2   &2      &0  &1   &1  &0   &0   &0    &0   &0               &0           \\
 1    &1   &1    &0   &0   &0   &2   &1      &1  &1    &1   &1 &0   &0    &0   &0                 &0      \\
 1    &3   &3    &3   &2   &1   &1   &2      &0  &0    &0   &0  &1  &0    &0   &0                &0       \\
 1    &2   &1    &1   &0   &0   &1   &2      &1  &0    &0   &0  &1   &1  &0   &0                &0         \\
 1    &1   &0    &0   &0   &0   &1   &1      &1  &1    &0   &0   &0   &1  &1    &0                 &0     \\
 1    &0   &0    &0   &0   &0   &1   &0      &1  &1    &0   &1   &0   &0  &1   &1        &0        \\
 1    &1   &1    &0   &1   &0   &1   &1      &0  &2    &0    &0   &0   &0  &1   &0     &1  \\
\endmatrix\right)$$
with rows indexed from left to right and columns indexed from up to down
by the elements of $M_0(S_5)$  in the order
$$(1,1),(1,\l^1),(1,\nu),(1,\l^2),(1,\nu'),(1,\l^3),(g_2,1),(g_2,r),(g_3,1),(g'_2,1),(g'_2,\e''),$$
$$(g_6,1),(g_2,\e),(g_3,\e),(g_4,1),(g_5,1),(g'_2,\e').$$
For $(x,\s)\in M_0(S_5)$, the coefficient of $(x',\s')\in M_0(S_5)$ in $\widehat{(x,\s)}\in\CC[M(S_5)]$  
is the entry of the matrix above in the row $(x,\s)$ and column $(x',\s')$;
the coefficient of any $(x',\s')\in M(S_5)-M_0(S_5)$ is $0$.

\subhead 3.9\endsubhead
The basis $\ti\BB_c$ defined above satisfies properties (I)-(V) in 0.1. (For (I) we use 3.1(g) and the results 
in 3.2.) It also satisfies the property stated in 0.2 (with the notion of primitive elements as in 3.2.)

\enddocument